\documentclass[11pt,reqno,tbtags]{amsart}
\usepackage{amssymb}
\usepackage{epsfig}
\usepackage{color}

\usepackage{colordvi}
\numberwithin{equation}{section}

\newtheorem{theorem}{Theorem}[section]
\newtheorem{lemma}[theorem]{Lemma}
\newtheorem{proposition}[theorem]{Proposition}
\newtheorem{corollary}[theorem]{Corollary}
\newtheorem{conjecture}[theorem]{Conjecture}

\theoremstyle{definition}
\newtheorem{example}[theorem]{Example}

\newtheorem{remark}[theorem]{Remark}

\newtheorem*{ack}{Acknowledgement}

\theoremstyle{remark}

\newenvironment{romenumerate}{\begin{enumerate}
 }{\end{enumerate}}

\newcounter{thmenumerate}
\newenvironment{thmenumerate}
{\setcounter{thmenumerate}{0}%
 \def\item{\par
 \refstepcounter{thmenumerate}\textup{(\roman{thmenumerate})\enspace}}
}
{}

\newcounter{xenumerate}   

\newcommand\pfitem[1]{\par(#1):}

\newcommand\upref[1]{\textup{\ref{#1}}}
\newcommand{\refT}[1]{Theorem~\upref{#1}}
\newcommand{\refC}[1]{Corollary~\upref{#1}}
\newcommand{\refL}[1]{Lemma~\upref{#1}}
\newcommand{\refR}[1]{Remark~\upref{#1}}
\newcommand{\refS}[1]{Section~\upref{#1}}
\newcommand{\refSS}[1]{Subsection~\upref{#1}}
\newcommand{\refP}[1]{Proposition~\upref{#1}}

\newcommand{\refConj}[1]{Conjecture~\upref{#1}}
\newcommand{\refand}[2]{\upref{#1} and~\upref{#2}}

\begingroup
  \divide\count255 by 60
  \multiply\count255 by -60
  \advance\count255 by \time
  \ifnum \count255 < 10 \xdef\klockan{\the\count1.0\the\count255}
  \else\xdef\klockan{\the\count1.\the\count255}\fi
\endgroup

\newcommand\nopf{\qed}   
\newcommand\noqed{\renewcommand{\qed}{}} 

\newcounter{CC} 
\newcommand{\CC}{\stepcounter{CC}\CCx}
\newcommand{\CCx}{C_{\arabic{CC}}}
\newcounter{cc}
\newcommand{\cc}{\stepcounter{cc}\ccx}
\newcommand{\ccx}{c_{\arabic{cc}}}
\newcounter{Ca}
\newcommand{\Ca}{\stepcounter{Ca}\Cax}
\newcommand{\Cax}{C_{\arabic{Ca}}(a)}
\newcommand{\Caq}{C(a)}
\newcounter{Cga}
\newcommand{\Cga}{\stepcounter{Cga}\Cgax}
\newcommand{\Cgax}{C_{\arabic{Cga}}(\ga)}
\newcommand{\Cgaq}{C(\ga)}

\newcommand\set[1]{\ensuremath{\{#1\}}}

\newcommand\bigpar[1]{\bigl(#1\bigr)}
\newcommand\Bigpar[1]{\Bigl(#1\Bigr)}
\newcommand\biggpar[1]{\biggl(#1\biggr)}

\newcommand\bigsqpar[1]{\bigl[#1\bigr]}

\def\rompar(#1){\textup(#1\textup)}    

\def\xexp(#1){e^{#1}}

\newcommand\floor[1]{\lfloor#1\rfloor}

\newcommand\ntoo{\ensuremath{{n\to\infty}}}

\newcommand\norm[1]{\|#1\|}
\newcommand\normm[2]{\|#1\|_{#2}}

\newcommand\iid{i.i.d.\spacefactor=1000}    
\newcommand\ie{i.e.\spacefactor=1000}
\newcommand\eg{e.g.\spacefactor=1000}
\newcommand\viz{viz.\spacefactor=1000}
\newcommand\cf{cf.\spacefactor=1000}
\newcommand{\as}{a.s.\spacefactor=1000}
\newcommand{\aex}{a.e.\spacefactor=1000}

\newcommand\ii{\mathrm{i}}

\newcommand{\tend}{\longrightarrow}
\newcommand\dto{\overset{\mathrm{d}}{\tend}}
\newcommand\pto{\overset{\mathrm{p}}{\tend}}

\newcommand\eqd{\overset{\mathrm{d}}{=}}

\newcommand\bbR{\mathbb R}
\newcommand\bbC{\mathbb C}

\newcommand\bbZ{\mathbb Z}
\newcommand{\bm}[1]{\mbox{\boldmath \ensuremath{#1}}}

\renewcommand\Re{\operatorname{Re}}
\renewcommand\Im{\operatorname{Im}}

\newcommand\E{\operatorname{\mathbb E{}}}
\renewcommand\P{\operatorname{\mathbb P{}}}
\newcommand\Var{\operatorname{Var}}
\newcommand\Cov{\operatorname{Cov}}

\newcommand\supp{\operatorname{supp}}

\newcommand\ga{\alpha}
\newcommand\gb{\beta}
\newcommand\gd{\delta}

\newcommand\gam{\gamma}
\newcommand\gG{\Gamma}

\newcommand\gl{\lambda}

\newcommand\gs{\sigma}
\newcommand\gss{\sigma^2}
\newcommand\eps{\varepsilon}

\newcommand\cD{\mathcal D}

\newcommand\cS{{\mathcal S}}
\newcommand\cT{{\mathcal T}}

\newcommand\etta{\boldsymbol1} 
\def\[#1]{[\![#1]\!]}

\newcommand\qq{^{1/2}}
\newcommand\qqi{^{-1/2}}
\newcommand\qi{^{-1}}

\newcommand\qw{^{1/4}}
\newcommand\qwi{^{-1/4}}
\newcommand\qww{^{3/4}}
\newcommand\qwwi{^{-3/4}}

\renewcommand{\=}{:=}

\newcommand\intoi{\int_0^1}
\newcommand\intoo{\int_0^\infty}
\newcommand\intoooo{\int_{-\infty}^\infty}

\newcommand\dd{\,d}

\newcommand{\mgf}{moment generating function}

\newcommand\lhs{left-hand side}
\newcommand\rhs{right-hand side}

\newcommand{\stablex}[1]{$#1$-stable}

\newcommand\ftuu{F_2(t,e^{\ii u},e^{-\ii u})}
\newcommand\ftu{F_u(t)}
\newcommand\sqrtt{\sqrt{1-4t}}
\newcommand\cTn{\mathcal{T}_n}
\newcommand\cTx[1]{\cT^{#1}}
\newcommand{\GW}{Galton--Watson}
\newcommand{\GWp}{\GW{} process}

\newcommand{\cGWt}{conditioned \GW{} tree}
\newcommand\ISE{\textup{ISE}}
\newcommand\tn{T_n}

\newcommand\bex{B^{\mathrm{ex}}}
\newcommand\ise{_{\text{\sc ise}}} 
\newcommand\muise{\mu\ise}
\newcommand\hmuise{\widehat{\mu}\ise}
\newcommand\fise{f\ise}
\newcommand\hfise{\widehat{f}\ise}

\newcommand\nuexc{\nu_{\text{\sc exc}}} 

\newcommand\tW{\overline{W}}
\newcommand\xW{\widetilde{W}}
\newcommand\tzeta{\tilde{\zeta}}
\newcommand\bx{\bar x}
\newcommand\by{\bar y}

\newcommand\bX{\bar X}
\newcommand\bXn{\bX_n}
\newcommand\hXn{\widehat{X_n}}

\newcommand\bY{\bar Y}
\newcommand\bYn{\bY_n}
\newcommand\hYn{\widehat{Y_n}}
\newcommand\hbYn{\widehat{\bYn}}
\newcommand\hgn{\widehat{g_n}}
\newcommand\hnu{\widehat{\nu}}
\newcommand\cor{C_0(\bbR)}
\newcommand\nun{\nu_{\gam n\qwi}}
\newcommand\hatf{\hat f}
\newcommand\CRT{continuum random tree}
\newcommand\BCRT{Brownian CRT}
\newcommand{\sob}[1]{\ensuremath{L^{2,#1}}}
\newcommand{\soba}{\sob{\ga}}
\newcommand\glk{\gl_1,\dots,\gl_k}
\newcommand\ssk{s_1,\dots,s_k}
\newcommand\dssk{\dd s_1 \dotsm \dd s_k}
\newcommand\xxk{x_1,\dots,x_k}
\newcommand\XXk{X_1,\dots,X_k}
\newcommand\yyk{y_1,\dots,y_k}
\newcommand\dyyk{\dd y_1 \dotsm \dd y_k}
\newcommand\gSzsk{\Sigma_{\zeta;\ssk}}
\newcommand\gSzxk{\Sigma_{\zeta;\xxk}}
\newcommand\gSzXk{\Sigma_{\zeta;\XXk}}
\newcommand\phizsk{\varphi_{\zeta;\ssk}}
\newcommand\phizxk{\varphi_{\zeta;\xxk}}
\newcommand\phizXk{\varphi_{\zeta;\XXk}}
\newcommand\hh{^{(h)}}
\newcommand\WW{V}
\newcommand\PP[1]{P_{#1}}
\newcommand\QQ[1]{Q_{#1}}

\newcommand{\Holder}{H\"older}
\newcommand{\Holdera}{\Holder($\ga$)}
\newcommand{\Holderb}{\Holder($\gb$)}
\newcommand\CS{Cauchy--Schwarz}

\newcommand\REM[1]{\texttt{[#1]}\marginal{XXX}}

\newcommand{\beq}{\begin{equation}}
\newcommand{\eeq}{\end{equation}}

\begin{document}
\title[density of the ISE and limit laws for embedded trees]
{The density of the ISE\\ and local limit laws for embedded trees}

\date{September 13, 2005} 

\author{Mireille Bousquet-M\'elou}
\address{CNRS, LaBRI, Universit\'e Bordeaux 1, 351 cours de la Lib\'eration,
  33405 Talence Cedex, France}
\email{mireille.bousquet@labri.fr}
\urladdr{http://www.labri.fr/Perso/\~{}bousquet/}

\author{Svante Janson}
\address{Department of Mathematics, Uppsala University, PO Box 480,
SE-751 06 Uppsala, Sweden}
\email{svante.janson@math.uu.se}
\urladdr{http://www.math.uu.se/\~{}svante/}


\begin{abstract} 
It has been known for a few years that the occupation measure of
several models of embedded trees converges, after a suitable
normalization, to the random measure called ISE (Integrated
SuperBrownian Excursion). 
Here, we prove  
a local version of this result: ISE has a 
(random)
H\"older continuous density,
and the {\em vertical profile\/} of embedded
trees converges to this density, at least for some such trees.

As a consequence, we derive a formula for the distribution of the
density of  ISE at a given point. This  follows from
earlier results 
by Bousquet-M\'elou 
on convergence of the vertical profile
at a fixed point.

We also provide a recurrence relation defining the moments of the
(random) moments of ISE.
\end{abstract}

\maketitle

\section{Introduction}\label{S:intro}

We consider some families of random labelled trees; the labels
will be integers (positive or negative).
Our main case is
binary trees, where each node is labelled 
with the difference between
the number of right steps and the number of left steps occurring in the
path from the root to the node. In particular, the root has label 0,
and the labels of two adjacent nodes differ by $\pm1$.
Note that the label of each node is simply its abscissa, if we draw
the tree in the plane in such a way that the right [left] child of a
node lies one unit to the right [left] of its parent.
We call this the \emph{natural labelling} of a binary tree.

\begin{figure}[ht]
\begin{center}
\input{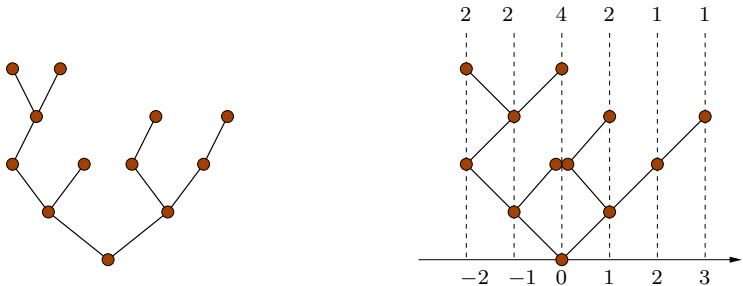}
\end{center}
\caption{A  binary tree having  vertical profile $[2,2,4; 2,1,1]$.}
\label{fig-profils}
\end{figure}

Given a labelled tree $T$, let $X(j;T)$ be the number of nodes in $T$ with
label $j$; the sequence $(X(j;T))_{j=-\infty}^{\infty}$ 
is the \emph{vertical profile} of the tree (Figure~\ref{fig-profils}).

Let $\tn$ be a random binary tree with $n$ nodes with the uniform
distribution, and let 
$X_n(j)\=X(j;\tn)$ be its vertical profile. 
It was shown by Marckert \cite[Theorem 5]{Marckert} that 
the (random) distribution of the labels in the tree converges, after
appropriate normalization, to the ISE (integrated superbrownian excursion)
introduced by Aldous \cite{AldousISE}, see also \cite{BM}.
The ISE is a
random probability measure; to emphasize this we will usually
write it as $\muise$.
(Actually, the result in \cite{Marckert} is stated for complete binary
trees, 
\ie{} binary trees where each node has either 0 or 2 children,
but the result transfers immediately by considering internal
nodes only; see at the end of Section~\ref{Sgf} for details.)
Marckert's result can be stated as follows, where 
$\gam\=2^{-1/4}$,
$\ell(v)$ denotes the label of $v$, 
and $\gd_x$ is the Dirac measure at $x$, 
\begin{equation}\label{sofie}
  \frac1n\sum_{j=-\infty}^\infty X_n(j) \gd_{\gam n\qwi j}
=
  \frac1n\sum_{v\in\tn} \gd_{\gam n\qwi \ell(v)}
\dto
\muise,
\qquad  \ntoo,
\end{equation}
with convergence in the space of probability measures on $\bbR$.
For complete binary trees, the result is the same, except that
now $\gam=1$ (and $n$ has to be odd).

Our first main result is a local version of \eqref{sofie}, 
showing that the vertical profile of random binary trees,
properly normalized, converges to the density $\fise$ of $\muise$; see
\refS{Sresults} for details.
Our second result consists of a recurrence relation
that characterizes the joint law of the moments of the ISE.

\begin{remark}
Different normalizations of $\muise$ are used in the literature. We
use the normalization of \cite{AldousISE}, also used by \eg{}
\cite{Marckert}.
The normalization in \cite{BM,MM} differs by a scale factor $2\qw$.
\end{remark}

Our local limit result actually holds for other families of labelled
(or: embedded) trees too.   Indeed,
the random measure $\muise$ 
arises naturally as a limit
for embedded trees in the following way \cite{AldousISE}.
Let $T_n$ be a random
\cGWt{} with $n$ nodes, 
\ie{} a random tree obtained as the family tree of a \GWp{}
conditioned on a given total population of $n$. 
(See \eg{} \cite{AldousII,Devroye} for details, and recall that this
includes \eg{} binary trees, complete binary trees,
plane trees and labelled (=Cayley) trees.
These random trees are also known as simply generated trees.)
The \GWp{} is defined using an offspring distribution;
let $\xi$ denote a random variable with this distribution.
We assume, as usual,
$\E \xi=1$ (the \GWp{} is critical)
and $ 0<\gss_\xi\=\Var\xi <\infty$.
Assign
\iid{} random 
variables $\eta_e$ to the edges of $\tn$.
We regard $\eta_e$ as the displacement from one endpoint of the edge
$e$ to the other, in the direction from parent to child; this gives
a labelling of the nodes such that the root has label 0 and each other
node $v$ has  label 
$\ell(v)\=\ell(v')+\eta_{vv'}$, where $v'$ is the parent of $v$.
For the purposes of this paper, we assume $\eta_e$ to be
integer valued. 
We further assume $\E\eta_e=0$ 
and $ 0<\gss_\eta\=\Var\eta_e <\infty$.
Define $X$ and $X_n$ as was done above for binary trees with their
natural labelling.
Then \eqref{sofie} holds, with $\gam\=\gs_\eta\qi\gs_\xi\qq$
\cite{AldousISE}, see also \cite{SJ164}.

We conjecture that a local version of \eqref{sofie} holds in this
generality, provided $\eta_e$ is not supported on a subgroup $d\bbZ$
of the integers with $d\ge2$,
but we will only prove this for two special cases, \viz{}
random plane trees with $\eta_e$ uniformly distributed on either
\set{\pm1} or \set{-1,0,1}, see \refT{Tplane}.

We state in \refS{Sise} some properties of the (random) density
function $\fise$, in particular that it exists. The proofs are given
in \refS{Spfise} after some preliminaries on the Brownian snake and
the \BCRT{} 
(\CRT) in \refS{Sbrown}.
Our results on the local limit law are stated in \refS{Sresults} 
and proved in Sections 
\ref{Spf1}--\ref{Strees}. 
Some further computations of (mixed) moments of the density $\fise(\gl)$
are given in \refS{Smom}. Our results on  moments of $\muise$ are
stated in \refS{Smom2} and proved  in \refS{Smom3}.

\smallskip

All unspecified limits below are as \ntoo.
We will use $C$ and $c$ with various subscripts
to denote various positive constants, not
depending on $n$ or other variables; for constants depending on a
parameter we use $C(a)$ and so on.

\begin{ack}
This research was mainly done at the Mittag-Leffler Institute,
Djursholm, Sweden, during the semester on Algebraic Combinatorics.
MBM  was partially supported by the European Commission's IHRP
  Programme, grant HPRN-CT-2001-00272, ``Algebraic Combinatorics in
  Europe''.

We further thank 
Ingemar Kaj and
Jean-Fran\c{c}ois Le Gall for helpful comments.
\end{ack}

\section{The density of the ISE}\label{Sise}

It is no surprise that the random measure $\muise$ is absolutely
continuous; the following theorem may well be known to experts, but we
have not found an explicit reference. 
(Related results for super-Brownian motion have been given by
\cite{KonnoShiga,Reimers,Sugitani}. It seems to be possible
but non-trivial to derive the existence of a density for ISE from
these results.)
We give a proof in \refS{Spfise}.

\begin{theorem}\label{T1}
\ISE{} has a \Holder{} continuous density.
In other words, there exists a continuous stochastic process $\fise(x)$,
$-\infty<x<\infty$, such that $\dd\muise(x)=\fise(x)\dd x$.
Moreover, the random function $\fise(x)$ has \as{} the following
properties:
\begin{romenumerate}
\item
$\fise$ has compact support: $\sup\set{|x|:\fise(x)\neq0}<\infty$;
\item
$\fise$ is \Holdera-continuous for every $\ga<1$;
\item
$\fise$ has a derivative $\fise'(x)$ \aex{} 
and in distribution sense,
and $\fise'\in L^p(\dd x)$ for every $p$ with $2\le p<\infty$.
\end{romenumerate}
\end{theorem}

Of course, the support of $\fise$ is random; (i) says that
there exists a random $M<\infty$ such that  $\fise(x)=0$ for $|x|>M$,
but no deterministic $M$ will do.

\begin{remark}\label{Rfise}
More precisely, the proof in \refS{Spfise}
shows that $\fise$ belongs to the generalized
Sobolev space $\soba$ for any $\ga<3/2$.
Loosely speaking, $\fise$ thus has ``$\ga$~derivatives in $L^2$''
for every real $\ga<3/2$. 
\end{remark}

Parts (ii) and (iii) of \refT{T1} come close to showing that $\fise$
has a continuous derivative, but we have not been able to prove it.
Indeed, it seems likely that the (fractional) derivatives in $L^2$
asserted by \refR{Rfise} are continuous. Hence we make the following
conjecture. 
\begin{conjecture}
The density $\fise$ has \as{} a continuous derivative, but not a
second derivative. 
\end{conjecture}

The marginal distributions of $\fise$, \ie{} the distributions of
$\fise(\gl)$ for fixed $\gl$, will be described in Corollaries
\refand{Cfise}{Cfise0}. 
Moments and mixed moments of $\fise(\gl)$ will be computed in
\refS{Smom}.

\section{Local limit results}\label{Sresults}

Our main result is the following local limit result for naturally
embedded random binary trees, conjectured in
\cite{BM}.

We let $\bXn(x)$ denote the function obtained by extending $X_n(j)$
to arbitrary real arguments by linear interpolation; thus
$\bXn(j)=X_n(j)$ for every integer $j$, and $\bXn$ is linear on each
interval $[j,j+1]$.

$\cor$ denotes, as usual, the Banach space of continuous functions on
$\bbR$ that tend to 0 at $\pm\infty$. We equip $\cor$ with the usual
uniform topology defined by the supremum norm.

Recall that we have defined the constant $\gamma$ as $2\qwi$ for
binary trees and~$1$ for complete binary trees.

\begin{theorem}\label{T2}
Consider random binary trees or random complete binary trees with their
natural labelling. Then, 
as \ntoo,
  \begin{equation}\label{t2a}
	\frac1n \gam\qi n\qw \bXn\bigpar{\gam\qi n\qw x}
\dto \fise(x),
  \end{equation}
in the space $\cor$ with the usual uniform topology.
Equivalently,
\begin{equation}\label{t2b}
 n\qwwi \bXn\bigpar{n\qw x}
\dto \gam\fise(\gam x).
\end{equation}
\end{theorem}

Note that the functions on the \lhs{s} of \eqref{t2a} and \eqref{t2b}
are density functions, \ie{}
non-negative functions with integral 1.
Proofs will be given in Sections \ref{Spf1}--\ref{Slastpf}.

\begin{corollary}\label{C2}
For random binary trees or random complete binary trees with their
natural labelling,
if \ntoo{} and $j_n/n\qw \to x$, where
$-\infty<x<\infty$,
then 
$n\qwwi X_n(j_n)\dto \gam \fise(\gam x)$.
\end{corollary}

It follows by combining this with results in 
Bousquet-M\'elou \cite{BM}
that the marginal distributions of $\fise$ are as
conjectured there.

\begin{corollary}\label{Cfise}
For every real $x$, the distribution of $\fise(x)$ 
is given by the \mgf{}
\begin{equation*}
\E e^{a \fise(x)} = L(2\qwi |x|,2\qwi a), 
\qquad |a|<2^{2+1/4} 3\qqi,  
\end{equation*}
where
$$
L(x, a)\=
1+ \frac {48}{\ii\sqrt \pi }\int_{\Gamma}
\frac{A(a/v^3)e^{-2x v}}
{(1+A(a/v^3)e^{-2x v})^2}v^{5} e^{v^4} dv,
\qquad x\ge0,
$$
$A(y)\equiv A$ is the unique solution of
\begin{equation*}
A=\frac y{24} \frac{(1+A)^3}{1-A}
\end{equation*}
satisfying $A(0)=0$,
and the integral is taken over 
$$
\Gamma= \set{1-te^{-\ii\pi/4}, t\in (-\infty,0]} \cup 
\set{1+te^{\ii\pi/4},t\in [0,\infty)}. 
$$
\end{corollary}

In particular, the density at $x=0$  has a simple law. (See again \cite{BM}.)

\begin{corollary}\label{Cfise0}
  $\fise(0)$ has the same distribution as $2\qw 3\qi T\qqi$, 
where $T$ is a
positive \stablex{2/3} variable with Laplace transform 
$\E e^{-t T}=e^{-t^{2/3}}$.

Hence $\fise(0)$ has the moments 
\begin{equation*}
 \E \fise(0)^r=2^{r/4}3^{-r} \frac{\Gamma(3r/4+1)}{\Gamma(r/2+1)},
\qquad -4/3<r<\infty.
\end{equation*}
\end{corollary}

As said in the introduction,
we conjecture that the local limit results hold also
for \cGWt{s} with random
labellings defined by \iid{} random increments $\eta_e$ along the
edges;
a precise formulation is as follows.
Recall that the \emph{span} of $\eta_e$ is the largest integer $d\ge1$
such that $\eta_e$ \as{} is a multiple of $d$.

\begin{conjecture}\label{conj1}
Consider a random \cGWt{} $\tn$
with a random labelling defined as above by integer
valued random variables $\eta_e$ with mean $0$, finite variance
$\gss_\eta>0$ and span $1$. 
Then, the conclusions \eqref{t2a} and \eqref{t2b} of \refT{T2} hold,
with $\gam\=\gs_\eta\qi\gs_\xi\qq$.
\end{conjecture}

If this conjecture holds, the conclusion of \refC{C2} holds too.

As said in the introduction, we can prove the conjecture in two
special cases, both considered in \cite{BM}.

\begin{theorem}
  \label{Tplane}
\refConj{conj1} holds if\/ $\tn$ is a random plane tree and
$\eta_e$ is uniformly chosen at random from \set{\pm1} or from
\set{-1,0,1}.
\end{theorem}
For these two cases, $\gss_\xi=2$ and $\Var\eta_e=1$ and $2/3$;
hence $\gam=2\qw$ and $\gam=2\qwi3\qq$, respectively.

\begin{remark}
It follows from the proof of \refT{T2} in \refS{Spf1} that 
to prove \refConj{conj1} in further cases,
it suffices to prove the estimate in \refL{L3}.
\end{remark}

\section{The moments of  ISE}
\label{Smom2}
Let $T_n$ be a random binary tree with $n$ nodes, and let $\mu_n$ be
the following (random) probability distribution:
\begin{equation}
\mu_n=
  \frac1n\sum_{v\in\tn} \gd_{(2 n)\qwi \ell(v)}.
\end{equation}
As recalled in the introduction, $\mu_n$ converges to $\muise$. 
The $i$th moment of $\mu_n$, denoted $m_{i,n}$, is itself a random
variable:
$$
m_{i,n}=   \frac1n\sum_{v\in\tn} (2n)^{-i/4}\ell(v)^i= 2^{-i/4}
n^{-1-i/4} \sum_{v\in\tn}\ell(v)^i.
$$
We shall prove that the sequence  $m_{1,n}, m_{2,n},
\ldots$ converges in distribution to the sequence $m_1, m_2, \ldots$
of moments of ISE, and  compute the joint moments of the $m_i$: 
$$
\E( m_{1}^{p_1} m_{2}^{p_2} \cdots m_{r}^{p_r} ),
$$
for all (fixed) values of $p_1, p_2, \ldots , p_r$. 
The moments of 
the  $m_{i}$, being the moments of the moments (of $\muise$) should
probably be called the {\em grand-moments\/} of $\muise$. 
Note that the grand-moments of a random probability measure, provided
they do not grow too quickly, determine the distribution of the
sequence of moments of the measure, and thus the distribution of the
random measure. 

In order to state our result, we introduce some
notation. A {\em  partition\/} $\lambda$ 
of an integer $k$ is a sequence $(\lambda_1, \ldots , \lambda_p)$ of
non-decreasing positive integers summing to $k$. The value $k$ is
called the \emph{weight} of $\lambda$, also denoted $k=|\lambda|$. 
For instance,
$\lambda=(1,1,3,4) $ is a partition of $k=9$. The $\lambda_i$ are
called the parts of $\lambda$. We shall also use {\em
  extended partitions\/}, in which the positivity condition on the
parts is  relaxed by simply requiring that $\lambda_i$ is
non-negative. Hence $\lambda=(0,0,1,1,3,4) $ is an extended partition
of $9$. The {\em union\/} $\sigma \cup \tau$ of two extended partitions
$\sigma=(\sigma_1, \ldots , \sigma_p)$ and
$\tau=(\tau_1, \ldots , \tau_q)$ is obtained by reordering the sequence
$(\sigma_1, \ldots , \sigma_p,\tau_1, \ldots , \tau_q)$. For any
$p$-tuple $(\sigma_1, \ldots, \sigma_p)$ of non-negative integers, we
denote by $\bar \sigma$ the extended partition obtained by reordering
the $\sigma_i$.
Given two
$p$-tuples $\sigma$ and $\lambda$, 
we write $\sigma\le \lambda$ if $\sigma_i \le \lambda_i$ for
 $1\le i\le p$. 

We shall denote
\beq
\label{moments-notation}
m_{1,n}^{p_1} m_{2,n}^{p_2} \cdots m_{r,n}^{p_r} := m_{\lambda,n}
\quad \hbox{and} \quad 
m_{1}^{p_1} m_{2}^{p_2} \cdots m_{r}^{p_r} := m_{\lambda} 
\eeq
where $\lambda=1^{p_1}2^{p_2}\cdots$ 
 is the partition having $p_1$ parts equal to $1$,
$p_2$ parts equal to $2$ and so on. 
The value of $\E(m_{\lambda})$ will be expressed in terms of a
   rational number   $c_\lambda$, which we actually define
   for   any \emph{extended} partition $\lambda=(\lambda_1, \ldots ,
   \lambda_p)$. The definition works  by induction on   $p+|\lambda|$
   as follows:
\begin{itemize}
\item $c_\emptyset = -2$,
\item  $c_\lambda=0\ $ if $|\lambda|$ is odd,
\item $c_\lambda=  (p+|\lambda|/4-3/2) c_{\lambda'}\ $ if $
  \lambda_1=0$, with $\lambda'=(\lambda_2, \ldots , \lambda_p)$,
\item if $\lambda_1>0$, 
\beq
\label{c-lambda-rec}
c_\lambda= \frac 1 4 
\sum_{\emptyset \neq I \subsetneq [p]} 
c_{\lambda_I} c_{\lambda_J} 
+ \sum_{\sigma\le \lambda, |\sigma|=|\lambda|-2}
 \binom{\lambda}{\sigma}  c_{\bar\sigma }
\eeq
where $J=[p]\setminus I$, 
$\lambda _I= (\lambda_{i_1}, \lambda_{i_2}, \ldots , \lambda_{i_r})$ 
if $I=\{i_1, \ldots ,i_r\}$ with $1\le i_1 < \cdots < i_r \le p$, 
the second sum runs over
all non-negative $p$-tuples $\sigma$ (not necessarily partitions)
satisfying the two required conditions, and 
$\binom{\lambda}{\sigma} = \prod_{i=1}^p  \binom{\lambda_i}{\sigma_i} $.
\end{itemize}

\begin{theorem}
\label{thm-grandmoments} 
As $n\rightarrow \infty$, the moments
$m_{1,n}, m_{2,n}, \ldots$ of the occupation measure of binary trees
converge jointly 
in distribution to the moments $m_1, m_2, \ldots$ of ISE. The
convergence of moments holds as well, and for all partitions $\lambda$,
the joint
$\lambda$-moment of the random variables $m_{i,n}$, defined
by~\eqref{moments-notation},  satisfies 
$$
\E(m_{\lambda,n}) =\E(m_{\lambda}) =0 \quad \hbox{ if\/ } |\lambda|
\text{ is odd},
$$
and otherwise
\begin{equation}
  \label{t41}
\E(m_{\lambda,n}) \rightarrow \E(m_{\lambda})=\frac{2^{-|\lambda|/4}\,
c_\lambda\,
  \Gamma(1/2)}{ \Gamma(p+|\lambda|/4-1/2)},
\end{equation}
where the number $c_\lambda$ is defined just above.
\end{theorem}
The vanishing  of $\E(m_{\lambda,n})$ when $ |\lambda|$ is odd, is a
straightforward consequence of the symmetry of $T_n$. 
The proof is given in Section~\ref{Smom3}.

\begin{example}[{\bf The average moments of ISE}]
When $\lambda$ has a single part, equal to $2k$ for $k\ge 1$, the
above recurrence relation  gives $c_{(2k)}= k(2k-1)c_{(2k-2)}$,
together with the
  initial condition $c_{(0)}= (1-3/2)c_\emptyset=1$. Hence the mean
  of the $2k$th moment of the random probability $\mu_n$ satisfies
\beq\label{2kth}
\E(m_{2k,n}) \rightarrow \E(m_{2k}) =\frac{(2k)! \, \Gamma(1/2)}{2^{3k/2} \,
  \Gamma((k+1)/2)}.
\eeq
\end{example}

\begin{example}[{\bf The moments of the average of ISE}]
\label{ex-first-moment}
Let $m_{1,n}$ denote the mean of $\mu_n$. Then $\E(m_{1,n}^p)=\E(m_{1}^p)=0$ if
$p$ is odd, while  
$$
\E\left(m_{1,n}^{2k}\right)\rightarrow \E\left(m_{1}^{2k}\right)
 =\frac{a_k \, \Gamma(1/2)}{2^{k/2} \,
  \Gamma((5k-1)/2)}, 
$$
where $a_0=-2$ and for $k>0$,
$$
4\, a_k=  \sum_{i=1}^{k-1} \binom{2k}{2i} a_i a_{k-i}
+ {k(2k-1)(5k-4)(5k-6)} a_{k-1}.
$$
Indeed, $a_k=c_\lambda$, where 
$\lambda=1^{2k}$,  and the recurrence  relation~\eqref{c-lambda-rec}
translates into the above recursive definition of $a_k$. Note that
each $2k$-tuple $\sigma$ occurring in~\eqref{c-lambda-rec} contains 2
coefficients equal to zero, so that we also use the part of the
definition of $c_\lambda$ that deals with the case $\lambda_1=0$. This
value of 
$\E(m_{1}^{2k})$ was already obtained in~\cite{SJ176}. 
 \end{example}

\begin{example}[{\bf The first two moments of ISE}] Let us finally
  work out the joint distribution of $m_1$ and $m_2$. We have
  $\E(m_{1,n}^{2k+1} m_{2,n}^\ell) = \E(m_{1}^{2k+1}
  m_{2}^\ell) =0$, and 
$$
\E(m_{1,n}^{2k} m_{2,n}^\ell) \rightarrow \E(m_{1}^{2k}
  m_{2}^\ell) = 
  \frac{a_{k,\ell} \, \Gamma(1/2)}{2^{(k+\ell)/2} \,
  \Gamma((5k+3\ell -1)/2)}, 
$$
where $a_{0,0}=-2$ and the $a_{k,\ell}$ are determined by induction
on $k+\ell$: 
\begin{multline*}
  a_{k,\ell} = \frac 1 4 \sum_{(0,0) < (i,j) <(k,\ell)} \binom{2k}{2i}
\binom{\ell}{j} a_{i,j} a_{k-i,\ell-j}
+
2\ell(\ell-1) a_{k+1,\ell-2}
\\+
\frac 1 4 k (2k-1)(5k+3\ell-4)(5k+3\ell-6)a_{k-1,\ell}
+
\frac 1 2   (4k+1)\ell(5k+3\ell-4)a_{k,\ell-1}.
\end{multline*}
Here, $a_{k,\ell}=c_\lambda$ with $\lambda=1^{2k}2^\ell$. 
In the right-hand side of the equation, the second (resp.\ third,
fourth) term corresponds to the case $\bar \sigma=  1^{2k+2}
2^{\ell-2}$ (resp.\ $\bar \sigma=0^21^{2k-2} 2^\ell$, $\bar
\sigma=01^{2k} 2^{\ell-1}$). The last case occurs both when we
replace a part of $\lambda$ equal to 2 by a zero part, and when we
decrease by 1  a part equal to 1 \emph{and} a part equal to 2. 
Of course, this generalizes Example~\ref{ex-first-moment} (which
corresponds to $\ell=0$). 
\end{example}

It seems likely that \refT{thm-grandmoments} extends to randomly 
labelled \cGWt{s} as in \refConj{conj1}, at least under some moment
conditions on $\xi$ and $\eta_e$, 
where the measure $\mu_n$ is defined by the left-hand side 
of~\eqref{sofie} and, as usual,  $\gam\=\gs_\eta\qi\gs_\xi\qq$.
We show this for the
special case in \refT{Tplane}. 

\begin{theorem}
\label{thm-grandmoments-plane} 
If\/ $\tn$ is a random plane tree and
$\eta_e$ is uniformly chosen at random from \set{\pm1} or from
\set{-1,0,1},
then the conclusions of
\refT{thm-grandmoments} hold, 
where the measure $\mu_n$ is defined by~{\rm{\eqref{sofie}}}, with
$\gamma=2^{1/4}$ and $\gamma=2^{-1/4}3^{1/2}$, respectively.
\end{theorem}

To conclude this section, we want to underline briefly a similarity
between the density of ISE and the local time of the (normalized) Brownian
excursion.  
In fact, \refT{T1} shows that the vertical profile of a random binary
tree converges, after suitable rescaling, to the density of the ISE.
Similarly, as shown by Drmota and Gittenberger \cite{DrmotaG},
the horizontal profile converges to the local time of the Brownian excursion.

We can develop this analogy to grand-moments as follows.
Consider again the random binary tree $T_n$.
Let $d(v)$ denote the depth (the distance from the root) of vertex $v$
and define the probability measure
\begin{equation}
  \label{nun}
\nu_n := \frac 1 n \sum_{v \in T_n} \delta_{2^{-3/2}n^{-1/2} d(v)},
\end{equation}
describing the horizontal profile, \ie{} the
distribution of the depths (after rescaling).
It is known, as an immediate consequence of Aldous \cite{AldousII,AldousIII},
that as \ntoo, $\nu_n\dto \nuexc$, where the random probability measure
$\nuexc$ is the occupation measure of the Brownian excursion, and thus
has the local time of the Brownian excursion as density.

The similarity with the vertical profile and ISE 
is obvious, and we adopt below the same notation as before ($m_{i,n},
m_{\lambda,n}$, etc.)\ for the moments of $\nu_n$ and $\nuexc$. 
In particular, now
$m_i\=\int x^i \dd\nuexc(x)=\intoi e(t)^i\dd t$, where $e(t)$ is a
Brownian excursion.
Then a
result similar to Theorem~\ref{thm-grandmoments} holds:
\begin{theorem}
\label{thm-grandmoments-exc} 
As $n\rightarrow \infty$, the moments
$m_{1,n}, m_{2,n}, \ldots$ of the horizontal profile (depth
distribution) measure $\nu_n$ 
converge jointly in distribution to the moments $m_1, m_2, \ldots$ of
$\nuexc$. The 
convergence of moments holds as well, and for all partitions $\lambda$,
the joint
$\lambda$-moment of the random variables $m_{i,n}$, defined
by~\eqref{moments-notation},  satisfies 
\begin{equation}
  \label{t45}
\E(m_{\lambda,n}) \rightarrow \E(m_{\lambda})
 =\frac{2^{-3|\gl|/2}\, d_\lambda\,  \Gamma(1/2)}
  { \Gamma(p+|\lambda|/2-1/2)},
\end{equation}
where the number $d_\lambda$ is defined by
\begin{itemize}
\item $d_\emptyset = -2$,
\item $d_\lambda=  (p+|\lambda|/2-3/2) d_{\lambda'}\ $ if $
  \lambda_1=0$, with $\lambda'=(\lambda_2, \ldots , \lambda_p)$,
\item if $\lambda_1>0$, 
\beq
\label{d-lambda-rec}
d_\lambda= \frac 1 4 
\sum_{\emptyset \neq I  \subsetneq [p]} 
d_{\lambda_I} d_{\lambda_J} 
+  \sum_{\sigma\le \lambda,\, |\sigma|=|\lambda|-1} 
 \binom{\lambda}{\sigma} d_{\bar\sigma },
\eeq
\end{itemize}
with the same notation as in~{\rm\eqref{c-lambda-rec}}.
\end{theorem}
The proof is very similar to the proof of
Theorem~\ref{thm-grandmoments}, but simpler. Note that the binomial
coefficient $\binom{\lambda}{\sigma}$ is simply equal to one of the
$\lambda_i$. The proof is 
sketched at the end of Section~\ref{Smom3}.

The grand-moments in \eqref{t45} have been computed by 
a different method by Richard \cite{Richard}; a special case (moments
of $m_1$ and $m_2$) is given by Nguyen The \cite{Nguyen}, 
and the moments of the
Brownian excursion area $m_1$ were found already by Louchard \cite{louchard},
see also \cite{flajolet-louchard} and
\cite{SJ146}.
The grand-moments in \eqref{t45}, 
as well as the grand-moments in \eqref{t41} above,
can also be derived by the method of \cite[Section 5]{SJ146}, which is
related to the method used here but phrased in different terms.
(Presumably, the method of \cite{Richard} too applies to \eqref{t41}
as well.)

Again, the same result holds for random plane trees as well, provided 
we change the scale factor $2^{-3/2}$ in \eqref{nun} to $2^{-1/2}$.

\begin{remark}
\label{Rdyck}
A \emph{Dyck path} of length $2n$ is a 1-dimensional walk 
starting and ending at $0$, taking steps in $\{-1,+1\}$, and never
reaching a negative position.
There is a well-known correspondence between plane trees with $n+1$
vertices and Dyck paths of length $2n$, where the Dyck path gives the
depths of the vertices along
the depth-first walk on the tree.
It follows easily that
\refT{thm-grandmoments-exc} holds for moments
of a uniformly chosen random Dyck path $w_n$ of length $2n$ too defined by
$$
m_{k,n} := \frac 1{2n} \sum_{i=1}^{2n} (2n)^{-k/2} w_n(i)^k;
$$
this has previously been shown by Richard \cite{Richard}.
\end{remark}

\begin{remark} 
It is possible to use our methods to obtain results on grand-moments
of 
the vertical and horizontal occupation measures together, and thus 
on the joint distribution of the vertical and horizontal profile, and
also on the asymptotic distribution of the pair of labels
$\bigpar{\ell(v),d(v)}$. We leave this to the reader.
\end{remark}


\section{The Brownian snake and CRT}\label{Sbrown}

\subsection{The Brownian snake}\label{SSsnake}
We begin by recalling the definition of the Brownian snake, see 
Le Gall \cite[Chapter IV]{LeGall} or Le Gall and Weill \cite{LeGallW}
for further details; see also \cite[Section 4.1]{SJ176}.
Let $\zeta$, the \emph{lifetime}, be 
$\zeta\=2\bex$, where $\bex$ is a Brownian excursion on $[0,1]$. 
(In general, the lifetime $\zeta$ might be any (locally) \Holder{} continuous 
non-negative stochastic process on some interval $I$;
in other contexts, $\zeta$ is often taken to be reflected
Brownian motion on $[0,\infty)$ \cite{LeGall}.)
Let, for $s,t\in [0,1]$,
\begin{equation*}
m(s,t;\zeta) \= \min\set{\zeta(u):u\in[s,t]}
\qquad\text{when $s\le t$,} 
\end{equation*}
and $m(s,t;\zeta)\=m(t,s;\zeta)$ when $s>t$.
The Brownian snake with lifetime $\zeta$ then can be defined as the
continuous
stochastic process $W(s,t)$ on $[0,1]\times[0,\infty)$ such that,
conditioned on $\zeta$, $W$ is Gaussian with mean 0 and covariances
\begin{equation*}
  \Cov\bigpar{W(s_1,t_1)W(s_2,t_2)\mid\zeta}
=\min\bigpar{t_1,t_2,m(s_1,s_2;\zeta)}.
\end{equation*}
We have defined the Brownian snake as a random field with two
parameters, but we are really only interested in the specialization
$\tW(s)\=W(s,\zeta(s))$, $s\in[0,1]$;
this stochastic process is called the 
\emph{head} of the Brownian snake. 
(In fact, it is easily seen that the pair $(\zeta,\tW)$ determines
$W$; see further \cite{MM}.) 
Conditioned on $\zeta$, $\tW$ is a Gaussian process on $[0,1]$ with mean 0
and covariances 
$\E\bigpar{\tW(s)\tW(t)\mid\zeta}=m(s,t;\zeta)$.
Consequently, still conditioned on $\zeta$, $\tW(s)-\tW(t)$ 
has a normal distribution with mean 0 and variance
\begin{equation}
  \label{gss}
\Var\bigpar{\tW(s)-\tW(t)\mid\zeta}
=\gss(s,t;\zeta)\=\zeta(s)+\zeta(t)-    2     m(s,t;\zeta).
\end{equation}

The random probability measure $\muise$ 
can be defined
as the occupation measure of the process $\tW$, see \cite{LeGall},
\cite{LeGallW}
and the next subsection.
Hence, $\fise$ is the occupation density of $\tW$, also called its 
\emph{local time}.

\subsection{\BCRT}\label{SScrt}

The \BCRT{} (\CRT)
was introduced by Aldous \cite{AldousI,AldousII,AldousIII}
as a natural limit of rescaled finite random trees.
It is a random compact metric space that is a topological tree in the
sense that every pair of points $x,y$ are connected by a unique path
(homeomorphic to $[0,1]$), and that path has length $d(x,y)$. 
We let here and later 
$d$ denote the metric.
The \BCRT{} is further equipped with a probability measure $\nu$, which gives
a meaning to ``a random node'' in the CRT.

One of Aldous's characterizations of the \BCRT{}
\cite{AldousII,AldousIII,AldousISE} uses the 
distribution of the shape and edge lengths of the spanning
subtree $R_k$ spanned by the root $o$ and $k$ independent random nodes
$X_1,\dots,X_k$ in the tree. (Here $k$ is an arbitrary positive integer.)
Then \as, the subtree $R_k$ admits the root and $X_1,\dots,X_k$ as
leaves, and has exactly $k-1$ internal nodes, all of of degree 3; 
the leaves are
labelled but not the internal nodes.
If we ignore the edge lengths (which are positive real numbers), 
there are $(2k-3)!!$ possible
``shapes'' of $R_k$; for each shape we number the $2k-1$ edges in some
order. Letting $T^*_{2k}$ be the finite set of shapes,
$R_k$ can thus be described by a shape $\hat t\in T^*_{2k}$ and the
edge lengths $x_1,\dots,x_{2k-1}>0$, and 
for the \BCRT, $R_k$ has density
\cite[Lemma 21]{AldousIII}
\begin{equation}
  \label{rk}
f(\hat t;x_1,\dots,x_{2k-1})=s e^{-s^2/2},
\qquad s=\sum_{i=1}^{2k-1} x_i.
\end{equation}

Aldous \cite[Corollary 22]{AldousIII} 
also gives a construction of the \BCRT{} in terms of a (normalized)
Brownian excursion $\bex$. 
Let $\zeta=2\bex$. Then Aldous shows that there exists a 
function $\tzeta$ 
mapping $[0,1]$ onto the \BCRT, with the Lebesgue measure mapped to
$\nu$ and, \cf{} \eqref{gss},
\begin{equation}
  \label{tzeta}
d\bigpar{\tzeta(s),\tzeta(t)}
 =\zeta(s)+\zeta(t)-2m(s,t;\zeta)=\gss(s,t;\zeta).
\end{equation}
Indeed, the \BCRT{} can be defined as the quotient space of $[0,1]$
with the semi-metric $\gss(s,t;\zeta)$, identifying points of distance 0,
see \cite{LeGallW}.

The function $\tzeta$ is not injective, but if
$\tzeta(s)=\tzeta(t)$, and thus $\gss(s,t;\zeta)=0$, then
$\zeta(s)=\zeta(t)=m(s,t;\zeta)$, which implies $\tW(s)=\tW(t)$.
Hence we can define a continuous random function $\xW$ on the \BCRT{} by
$\xW(\tzeta(s))=\tW(s)$; conditioned on the CRT, the $\xW(x)$ are jointly
Gaussian with mean 0 and, by \eqref{gss} and \eqref{tzeta}, 
$\Var\bigpar{\xW(x)-\xW(y)}=d(x,y)$.
Thus $\xW$ is the random mapping of the \BCRT{} into $\bbR$ considered
by Aldous \cite{AldousISE}; Aldous defines ISE as the measure on
$\bbR$ that $\nu$ is mapped to by $\xW$. This is clearly the same as
the measure that $\tW$ maps Lebesgue measure on $[0,1]$ to, \ie{} the
occupation
measure of $\tW$ as claimed above.

\section{Existence of the density: Proof of \refT{T1}}\label{Spfise}

Although \refT{T1} follows easily from \refT{T2} and its proof, we
find it interesting to give a different, self-contained proof.
We use the standard Fourier method, see \eg{} \cite{occupation} and
the references there, together with Aldous's theory of the \BCRT{}
\cite{AldousII,AldousIII}. 
We define the Fourier transform $\widehat\mu$ of a finite measure $\mu$ by
$\widehat\mu(t)\=\int e^{\ii tx}\dd\mu(x)$.

\begin{lemma}\label{L1}
  If\/ $0\le \ga<3/2$, then
\begin{equation*}
\E\intoooo \bigpar{|t|^\ga |\hmuise(t)|}^2\dd t<\infty.
\end{equation*}
\end{lemma}

\begin{proof}
  Since $\muise$ is the occupation measure of $\tW$,
the head of the Brownian snake, its Fourier transform can be expressed
as
\begin{equation*}
  \hmuise(t)\=\intoooo e^{\ii t x}\dd\muise(x)
=\intoi e^{\ii t \tW(s)}\dd s.
\end{equation*}
Consequently, 
$|\hmuise(t)|^2 = \intoi\intoi e^{\ii t(\tW(s)- \tW(u))}\dd s\dd u$. 
Conditioned on $\zeta$, $\tW(s)-\tW(u)$ is by \eqref{gss}
a Gaussian random variable with mean 0 and variance
$\gss(s,u;\zeta)$. Hence,
\begin{equation*}
  \begin{split}
\E\bigpar{|\hmuise(t)|^2\mid\zeta}
&=\intoi\intoi \E\Bigpar{e^{\ii t(\tW(u)- \tW(s))}\,\Big|\,\zeta}\dd s\dd u
\\&
=\intoi\intoi e^{-t^2\gss(s,u;\zeta)/2}\dd s\dd u
  \end{split}
\end{equation*}
and thus, letting $U_1$ and $U_2$ be independent uniform random
variables on $[0,1]$,
\begin{equation}\label{4c}
\E{|\hmuise(t)|^2}
=\E\intoi\intoi e^{-t^2\gss(s,u;\zeta)/2}\dd s\dd u
=\E e^{-t^2\gss(U_1,U_2;\zeta)/2}.
\end{equation}

Let $\tzeta$ be as in \refSS{SScrt}. Then $X_i\=\tzeta(U_i)$, $i=1,2$, 
are two
independent random nodes in the \BCRT, and \eqref{tzeta} shows that
\eqref{4c} can be written
\begin{equation}
\E|\hmuise(t)|^2 
=\E e^{-t^2 d(X_1,X_2)/2}
.
\end{equation}
For  $\ga\ge0$ we thus have, letting $y=t^2$,
\begin{equation}\label{4a}
  \begin{split}
\E\intoo\bigpar{t^\ga|\hmuise(t)|}^2 \dd t
&
=\E \intoo t^{2\ga}e^{-t^2 d(X_1,X_2)/2}\dd t
\\&
=\E \intoo \tfrac12 y^{\ga-1/2}e^{-y d(X_1,X_2)/2}\dd y
\\&
=\E \tfrac12  \bigpar{d(X_1,X_2)/2}^{-\ga-1/2} \Gamma(\ga+1/2)
\\&
=\Cga\E d(X_1,X_2)^{-\ga-1/2} .
  \end{split}
\end{equation}

{}From \eqref{rk} (with $k=2$) follows
the symmetry $d(X_1,X_2)\eqd d(X_1,o)$. Moreover, by the same formula
\eqref{rk} with $k=1$,
$d(X_1,o)$  has a Rayleigh distribution with density $x e^{-x^2/2}$.
Hence, 
\begin{equation*}
  \begin{split}
\E d(X_1,X_2)^{-\ga-1/2} 
=\E d(X_1,o)^{-\ga-1/2} 
=\intoo x^{-\ga-1/2} x e^{-x^2/2}\dd x
<\infty,
  \end{split}
\end{equation*}
 when $\ga<3/2$, and the result follows from \eqref{4a} and the
symmetry of $|\hmuise|$.
\end{proof}

By \refL{L1}, if $0\le\ga<3/2$, then
$\intoooo \bigpar{|y|^\ga |\hmuise(y)|}^2\dd y<\infty$ a.s. 
Taking first $\ga=0$, we see that $\hmuise\in L^2(\bbR)$; 
by Plancherel's theorem \cite[Theorem 7.9]{Rudin}
this shows that $\muise$ is absolutely continuous with a density
$\fise\in L^2$. Note that the Fourier transform $\hfise$ coincides
with $\hmuise$.

For $\ga\ge0$, we define the (generalized) Sobolev space \soba{}
by
\begin{equation}\label{sob}
  \soba \=\set{f\in L^2(\bbR):
\normm{f}{2,\ga}^2\=\intoooo \bigpar{(1+|t|)^\ga |\hatf(t)|}^2\dd t<\infty},
\end{equation}
where $\hatf$ is the Fourier transform of $f$.
\refL{L1} thus shows that \as{} $\fise\in\soba$ for every $\ga<3/2$.
(There is no problem with null sets, since it suffices to consider
rational $\ga$, say.)

Further, for $0<\ga<1$, we define
the \Holder{} space $H_\ga$
as the space of bounded continuous functions
$f$ on $\bbR$ such that
$|f(x)-f(y)|\le C |x-y|^\ga$ for some $C$ and all $x$ and $y$.

To show that $\fise$ is (\ie{} can be chosen) continuous with the
regularity properties  in \refT{T1}, we use some general embedding
properties of these spaces.

\begin{lemma}\label{Lsob}
  \begin{thmenumerate}
\item
  If\/ $0\le \ga<1/2$ and $1/2\ge1/p>1/2-\ga$, then $\soba\subset L^p$.
\item
  If\/ $1/2<\ga<3/2$, then $\soba\subset H_{\ga-1/2}$.
\item
  If\/ $\ga\ge1$ and $f\in\soba$, then $f$ has a derivative $f'$ in
  distribution sense and \aex, with $f'\in\sob{\ga-1}$.
  \end{thmenumerate}
\end{lemma}

This lemma is well-known: (i) and (ii) are special cases of the
Sobolev (or Besov) embedding 
theorem, see \eg{} \cite[Theorem 6.5.1]{BL} or \cite[Chapter V]{Stein}; 
indeed, we may also take
$1/p=1/2-\ga$ in (i). However, since the proof of the general
embedding theorem is quite technical, we give a simple proof of this
special case. 

\begin{proof}
  \pfitem{i}
We may assume $p>2$ since the case $p=2$ follows by Plancherel's theorem.
Define $p'\in(1,2]$ by $1/p'=1-1/p$. By \Holder's inequality,
\begin{equation*}
  \begin{split}
  \intoooo|\hatf|^{p'}
&\le 
\Bigpar{\intoooo\bigpar{(1+|t|)^{\ga p'}|\hatf(t)|^{p'}}^{2/p'}\dd t}^{p'/2}
\\&\hskip6em\times
\Bigpar{\intoooo\bigpar{(1+|t|)^{-\ga p'}}^{2/(2-p')}\dd t}^{1-p'/2}
\\&
= \norm{f}_{2,\ga}^{p'}
\Bigpar{\intoooo(1+|t|)^{-2\ga p'/(2-p')}\dd t}^{1-p'/2}	
<\infty,
  \end{split}
\end{equation*}
since it is easy to check that $2\ga p'>2-p'$ when $1/p>1/2-\ga$.
Consequently, $\hatf\in L^{p'}$, which by 
the Hausdorff--Young inequality yields $\hat{\hatf}\in L^p$.
By the inversion theorem for the Fourier transform (defined for
tempered distributions, say), this yields $f\in L^p$.

\pfitem{ii}
First, by \Holder's (\CS's) inequality,
\begin{equation*}
  \begin{split}
  \intoooo|\hatf|
\le 
\Bigpar{\intoooo\bigpar{(1+|t|)^{\ga}|\hatf(t)|}^{2}\dd t}^{1/2}
\Bigpar{\intoooo(1+|t|)^{-2\ga }\dd t}^{1/2}
<\infty,
  \end{split}
\end{equation*}
since $2\ga>1$.
Hence $f$ has an absolutely integrable Fourier transform, which shows
that $f$ is a continuous bounded function given by the inversion
formula
$f(x)=(2\pi)\qi\int e^{-\ii xt}\hatf(t)\dd t$.
Hence, for any $x$ and $h>0$,
using \Holder's inequality again,
\begin{equation*}
  \begin{split}
|f(x+h)-f(x)|
&
=
\frac1{2\pi}
\Bigl|\intoooo\bigpar{e^{-\ii(x+h)t}-e^{-\ii xt}}\hatf(t)\dd t\Bigr|
\\&
\le
\Bigpar{\intoooo\bigl|e^{\ii ht}-1\bigr|^2\,|t|^{-2\ga}\dd t}\qq
\Bigpar{\intoooo|t|^{2\ga}|\hatf(t)|^2\dd t}\qq
\\&
\le
\Bigpar{h^{2\ga-1}\intoooo\bigl|e^{\ii u}-1\bigr|^2\,|u|^{-2\ga}\dd u}\qq
\norm{f}_{2,\ga}
\\&
\le 
\Cga \norm{f}_{2,\ga}h^{\ga-1/2}.
  \end{split}
\end{equation*}

\pfitem{iii}
We have $\widehat{f'}(t)=-\ii t\hatf$, with $f'$ taken as a
distribution.
Since $f \in L^{2,\alpha}$, this shows that $\widehat{f'}\in L^2$ 
and thus $f'\in L^2$ by Plancherel's theorem.
Consequently, by elementary
distribution theory, the derivative exists
\aex{}, and equals the distributional derivative $f'$. 
Further, from the definition \eqref{sob}, $f'\in \sob{\ga-1}$.
\end{proof}

Since, as remarked above, \refL{L1} shows that \as{} $\fise\in\soba$ for
every $\ga<3/2$, \refT{T1}(ii) follows by \refL{Lsob}(ii), while
\refT{T1}(iii) follows by \refL{Lsob}(iii) and (i) (applied to
$\fise'$).

Finally, \refT{T1}(i) follows because $\muise$ has compact support,
\viz{} the image of the compact set $[0,1]$ 
by the continuous function $\tW$.
\qed

\section{Local limit law for the density: Proof of \refT{T2}}\label{Spf1}

The proof is based on the known convergence \eqref{sofie} of  random
measures.
To obtain the stronger result in \refT{T2}
on convergence of densities, we use 
a compactness argument as follows. We begin with a  measure-theoretic lemma.
Recall that a Polish space is a space with a topology that can be
defined by a complete separable metric.
For generalities on convergence of random elements of metric spaces
(equipped with their Borel $\gs$-fields),
see \eg{} Billingsley \cite{Billingsley} or Kallenberg \cite{Kallenberg}.
In particular, recall that a sequence $(W_n)$ of random variables in a
metric space $\cS$ is \emph{tight} if for every $\eps>0$, there exists
a compact subset $K\subseteq\cS$ such that $\P(W_n\in K)>1-\eps$ for
every $n$; in a Polish space, this is equivalent to relative
compactness (of the corresponding distributions)
by Prohorov's theorem \cite[Theorems 6.1 and 6.2]{Billingsley},
\cite[Theorem 16.3]{Kallenberg}.
Recall further that ``convergence in distribution'' really means
convergence of the corresponding distributions, but it is often convenient
to talk about random variables instead of their distributions.

\begin{lemma}
  \label{Lconv}
Let $\cS_1$ and $\cS_2$ be two Polish spaces, and
let
$\phi:\cS_1\to\cS_2$ be an injective continuous map.
If\/ $(W_n)$ is a tight sequence of random elements of $\cS_1$ such that
$\phi(W_n)\dto Z$ in $\cS_2$ for some random $Z\in\cS_2$, then
$W_n\dto W$ in $\cS_1$ for some $W$ with $\phi(W)\eqd Z$.
\end{lemma}

\begin{proof}
By Prohorov's theorem, 
each subsequence of $(W_n)$ has a subsequence that converges in
distribution to some limit.
Let $W'$ and $W''$ be limits in distribution
of two such subsequences $W_{n_i'}$ and $W_{n_i''}$. 
Since $\phi$ is continuous, $\phi(W_{n_i'})\dto \phi(W')$ and 
$\phi(W_{n_i''})\dto \phi(W'')$.
Hence, $\phi(W')\eqd Z \eqd \phi(W'')$.

Let $A$ be a (Borel) measurable subset of $\cS_1$.
By the Souslin--Lusin theorem 
\cite[Theorem III.21, see also  III.16--17]{DM},
$\phi(A)\subseteq\cS_2$ is measurable.
Thus, using the injectivity of $\phi$, 
$$
\P(W'\in A)=\P\bigpar {\phi(W')\in \phi(A)}
=\P\bigpar{\phi(W'')\in \phi(A)}=\P(W''\in A).
$$
Consequently, $W'\eqd W''$.

In other words, there is a unique distribution of the subsequence
limits.
Thus, if $W$ is one such limit, then 
every subsequence of $(W_n)$ has a subsequence that converges in
distribution to $W$; this is equivalent to $W_n\dto W$.
\end{proof}

Let $Y_n$ denote the random probability measure on the \lhs{} of
\eqref{sofie}, let $\nu_h$ be the probability measure with the
triangular density function $h\qi(1-|x|/h)_+$, and let $\bYn$ be the
convolution $Y_n*\nun$.
Note that $\bYn $ has the density 
$g_n(x)\= n\qi\gam\qi n\qw\bXn\bigpar{\gam\qi n\qw x} \in\cor$.
Since $Y_n\dto\muise$ 
by \eqref{sofie}, and $\nun\pto\gd_0$, it follows easily that
$\bYn\dto\muise$ too.

Let $\cS_1\=\set{f\in\cor:f\ge0}$, with the uniform topology inherited
from $\cor$,
let $\cS_2$ be the space of locally finite measures on $\bbR$
with the vague topology, 
see \eg{} Kallenberg \cite[Appendix A2]{Kallenberg},
and let $\phi$ map a function $f$ to the corresponding measure $f\dd x$,
\ie{}, $\phi(f)$ is the measure with density $f$. 
Then $\cS_1$ is a closed subset of the separable Banach space $\cor$, and
is thus Polish, and so is $\cS_2$ by \cite[Theorem A2.3]{Kallenberg}.
Further, $\phi$ is continuous and injective.

Take $W_n\=g_n$ in \refL{Lconv}. 
We have just shown that $\phi(g_n)=\bYn\dto\muise$ in 
the space of probability measures on $\bbR$
and thus also in the larger space $\cS_2$.
If we can show that the sequence $g_n$ is tight in $\cS_1$, or,
equivalently, in $\cor$, then \refL{Lconv} shows that $g_n\dto g$
for some random function $g\in\cor$, which further equals (in distribution) the
density $\fise$ of $\muise$; hence the conclusion of \refT{T2} follows.
It thus remains only to prove the following lemma.

\begin{lemma}
  \label{L2}
The sequence 
$g_n(x)\= n\qi\gam\qi n\qw\bXn\bigpar{\gam\qi n\qw x}$, $n=1,2,\dots$, 
is tight in $\cor$.
\end{lemma}

The central estimate in the proof of \refL{L2}, and thus of \refT{T2},
is the following, which will be proved in \refS{Sgf}.
For a sequence $x(j)$, we define its Fourier transform by 
$\widehat x(u)\=\sum_j x(j)e^{\ii ju}$; this equals the Fourier transform
of the measure $\sum_j x(j)\gd_j$ on $\bbR$.

\begin{lemma}
  \label{L3}
There exists a constant $\CC$ such that for all $n\ge1$ and $u\in[-\pi,\pi]$,
\begin{equation} \label{l3}
\E|n\qi\hXn (u)|^2 \le \frac{\CCx}{1+n u^4}.  
\end{equation}
\end{lemma}

We can now prove \refL{L2} as follows.
We have $\hYn(y)=n\qi \hXn\bigpar{\gam n\qwi y}$.
Consequently, $\hYn$ is a periodic function with period
$2\pi\gam\qi n\qw $,
and \refL{L3} translates to
\begin{equation}\label{erika}
\E|\hYn (y)|^2 \le \frac{\CCx}{1+ \gamma^4 y^4},
\qquad |y|\le \gam\qi n\qw\pi.  
\end{equation}
Further,
\begin{equation}\label{magnus}
 \hgn(y)=\hbYn(y)=\hYn(y) \hnu_{\gamma n^{-1/4}}(y). 
\end{equation}

\begin{lemma}
  \label{L4}
Suppose that $0\le a<3$. Then there exists a constant $\Caq$ 
such that if $h>0$ and $f$ is a function with period
$2\pi/h$, then
\begin{equation*}
  \intoo |y|^a |f(y)|^2 |\hnu_h(y)|^2 \dd y
\le \Caq  \int_{-\pi/h}^{\pi/h} |y|^a |f(y)|^2  \dd y.
\end{equation*}
\end{lemma}

\begin{proof}
  By the change of variables $y\mapsto h^{-1}y$, 
we may assume $h=1$.
Then $\hnu_1(y)=(\sin(y/2)/(y/2))^2$. Hence, for $k\neq0$ and
$|y|\le\pi$,
\begin{equation*}
  |\hnu_1(y+2k\pi)|=\frac{\sin^2(y/2)}{(k\pi+y/2)^2}
\le\frac{y^2}{k^2}
\end{equation*}
and
\begin{equation*}
  \begin{split}
  \int_{(2k-1)\pi}^{(2k+1)\pi} |y|^a |f(y)|^2 |\hnu_1(y)|^2 \dd y
&\le
(3\pi |k|)^a   \int_{-\pi}^{\pi}  |f(y)|^2 \frac{y^4}{k^4} \dd y
\\&
\le 
\Ca  |k|^{a-4} \int_{-\pi}^{\pi} |y|^a |f(y)|^2  \dd y.	
  \end{split}
\end{equation*}
For the case $k=0$, we use instead the estimate 
$|\hnu_1(y)|\le1$.
The result follows by summing over all $k$.
\end{proof}

Let $h\=\gam n\qwi$. Then, by \refL{L4}, \eqref{magnus} and
\eqref{erika}, for any fixed $a$ with $0\le a<3$,
\begin{align*}
  \E \intoo |y|^a|\hgn(y)|^2\dd y
&\le
 \Caq \E \int_{-\pi/h}^{\pi/h} |y|^a |\hYn(y)|^2 \dd y
\\&
= 
\Caq  \int_{-\pi/h}^{\pi/h} |y|^a \E |\hYn(y)|^2 \dd y
\\&
\le 
\Ca\intoo \frac{|y|^a}{1+\gamma^4 y^4}\dd y
\le \Ca.
\end{align*}

We have proved the following, taking $a=2\ga$.
\begin{lemma}
  \label{L5}
If\/ $0\le \ga<3/2$, then\/
$\E\normm{ g_n}{2,\ga}^2\le \Cgaq$,
for some\/ $\Cgaq$ not depending on $n$.
\nopf
\end{lemma}

Next, fix $\beta\in(0,1)$, and let $\ga=\gb+1/2<3/2$.
For $A,M>0$, let $K_{M,A}$ be the set of all functions $f$ in
$\cor$ such that $f(x)=0$ for $|x|\ge M$ and $\normm{f}{2,\ga} \le A$.
By \refL{Lsob}(ii),
the functions in $K_{M,A}$ are
all \Holderb-continuous with uniformly bounded norm; they thus form an
equicontinuous family. We may regard $K_{M,A}$ as a subset of
$C[-M,M]$, the space of continuous functions on the compact interval
$[-M,M]$, and it follows by the Arzela--Ascoli theorem
\cite[A5]{Rudin}
that $K_{M,A}$
is a relatively compact subset of $C[-M,M]$, and thus of $\cor$ too.
(Note that the functions in $K_{M,A}$ all vanish at $\pm M$.)

Let $\eps>0$.
It follows from \refL{L5} that 
there exists $A$
such that $\P(\normm{ g_n}{2,\ga} >A)\le \eps/2$ for every $n$.
Moreover, Marckert \cite[Theorem 5]{Marckert} also showed that
\begin{equation}
  \label{jesper}
n\qwi \sup\set{|j|:X_n(j)\neq0}
=
n\qwi \sup\set{|\ell(v)|:v\in \tn}
\dto W
\end{equation}
for some random variable $W$.
It follows from \eqref{jesper} that there exists $M$
such that 
\begin{multline*}
 \P\bigpar{g_n(x)\neq 0\text{ for some $x$ with $|x|>M$}}
\\
=
 \P\bigpar{ X_n(j)\neq 0\text{ for some $j$ with $|j|>\gam\qi n\qw M-1$}}
<\eps/2. 	
\end{multline*}
Consequently, 
$\P(g_n\in K_{M,A})>1-\eps$ for every $n$, 
which shows that the sequence $(g_n)$
is tight.
This completes the proof of \refL{L2}, and thus of \refT{T2}, except
for the proof of \refL{L3}.

\begin{remark}
A more concrete alternative to the compactness argument (\refL{Lconv})
used above is to define regularizations of functions $f$ on $\bbR$ by
$f\hh(x)\=h\qi\int_x^{x+h}f$ for $h>0$. Note that
$\fise\hh(x)=h\qi\muise[x,x+h]$. For each fixed $h>0$, the 
version of \eqref{t2a}  with both left- and right-hand side regularized
holds since the corresponding distribution functions converge.
We may then let $h\to0$, using the \Holder{} estimate obtained by
Lemmas \refand{L5}{Lsob} together with \cite[Theorem 4.2]{Billingsley}. 
\end{remark}

\section{Proof of \refL{L3}}\label{Sgf}
It remains to prove \refL{L3}. We consider first the case
of binary trees.
We introduce the sequence of generating functions
\begin{equation}   \label{fk}
F_k(t,x_1,\dots,x_k)
 \=
 \sum_{T\in \cT} t^{|T|} \prod_{i=1}^k \Bigpar{\sum_{v\in T} x_i^{\ell(v)}},
\end{equation}
where $\cT$ is the family of all (possibly empty) binary trees and
$|T|$ is the number 
of nodes in $T$. Thus $F_k$ is a power series in $t$, with
coefficients in $\bbZ[x_1, \ldots, x_k, 1/x_1, \ldots, 1/x_k]$, the
ring of Laurent polynomials in the $x_i$ with integer coefficients.
For $k=0$, the product in 
the definition of $F_0$ reduces to $1$, so that $F_0$ is simply the
generating function of binary trees. In what follows, we often denote
$\bm x= (x_1, \ldots,
x_k)$ and $F_k(\bm x)= F_k(t,x_1,\dots,x_k)$. Moreover, for any subset
$I$ of $[k]:=\{1, 2, \ldots , k\}$, we denote
$\bm x _I= (x_{i_1}, x_{i_2}, \ldots , x_{i_p})$ if $I=\{i_1, \ldots ,
i_p\}$ with $i_1 < \cdots < i_p$.

\begin{proposition}
\label{propo-enum-binary}
 The series $F_k$ can be determined by induction on
  $k\ge 0$ using
 \begin{equation}
   \label{fkx}
F_k( \bm x)= \etta_{[k=0]} + t \sum_{(I,J)} \left(\prod_{i \in I} \bx_i\right)
\left(\prod_{j\in J} x_j\right) F_{|I|} (x_I)F_{|J|} (x_J)
 \end{equation}
where the sum runs over all ordered pairs $(I,J)$ of subsets of $[k]$
such that $I \cap J =\emptyset$, and $\bx_i=1/x_i$.
In particular, 
\beq
\label{F0}
F_0= \frac{1-\sqrt{1-4t}}{2t}
\eeq
and each $F_k( \bm x)$ admits a rational expression in $F_0$ and the $x_i$.
\end{proposition}
\begin{proof}
The equation satisfied by $F_0$ reads $F_0=1+tF_0^2$ and is of course
very classical: it is 
obtained by splitting a binary tree into its left and right subtrees.
Note that the empty binary tree does not contribute to $F_k$ when $k>0$. Then,
every non-empty binary tree is formed of a root with  a left
subtree $T_1$ and a right subtree $T_2$. Hence, for $k\ge 1$,
\begin{align*}
F_k(\bm x) &= \sum_{T_1, T_2} t^{1+|T_1|+|T_2|} 
\prod_{i=1}^k \left( 1 + \sum_{v\in T_1} x_i^{\ell(v)-1}+ \sum_{v\in
  T_2} x_i^{\ell(v)+1} \right)\\
&=\sum_{T_1, T_2} t^{1+|T_1|+|T_2|}  \sum_{(I,J)}\prod_{i \in I}
 \left(  \sum_{v\in T_1} x_i^{\ell(v)-1}\right)
\prod_{j \in J}
 \left(  \sum_{v\in T_2} x_j^{\ell(v)+1}\right)
\end{align*}
where the sets $I$ and $J$ are as in the statement of the proposition.
The result follows upon exchanging the two sums.
\end{proof}
Actually, for the proof of \refL{L3} we need only a special case of
$F_2$. The above proposition gives  a simple explicit expression of
$F_2(t,x,y)$ in terms of $F_0$ (and $x$ and $y$), or, equivalently,  in terms
of the generating function  $B=F_0-1$ of non-empty binary trees: 
\begin{equation}
  \label{B}
B=
B(t)
=\frac{1-2t-\sqrt{1-4t}}{2t}.
\end{equation}

\begin{corollary}
  \label{CF2}
For any real $u$,
\begin{equation}\label{lf}
\ftuu
=
\frac{B(1+B)(1+2B-B^2)}{(1-B)(1+B-2B \cos u)^2} .
  \end{equation}
\end{corollary}

\begin{proof}
 The cases $k=1$ and $k=2$ of the previous proposition give
$$
F_1(x)=t F_0^2+ t(x+\bx) F_0 F_1(x)$$
and 
\begin{multline*}
F_2(x,y)= 
t F_0^2+ t(x+\bx)F_0F_1(x) + t(y+\by) F_0F_1(y) 
\\
+ t(\bx y + x \by) F_1(x) F_1(y) + t(xy+\bx \by ) F_0F_2(x,y).  
\end{multline*}
Using $F_0=1+B$ and $t=B/(1+B)^2$, this gives
\begin{gather*}
F_1(x)= \frac{B(1+B)}{1+B(1-x-\bx)},
\\
F_2(x,y)=\frac{B(1+B) (1+2B+B^2(1-xy-\bx\by))}
{(1+B(1-x-\bx))(1+B(1-y-\by))(1+B(1-xy-\bx\by))}.    
\end{gather*}
Specializing to $x=1/y=e^{\ii u}$ provides the result.
\end{proof}

By definition, 
\begin{equation}\label{emma}
\hXn(u)
=\sum_j X(j;T_n) e^{\ii j u}
=\sum_{v\in T_n} e^{\ii \ell(v) u}.
\end{equation}
Hence, if $\cTn\=\set{T\in\cT:|T|=n}$ 
is the family of binary trees of size $n$,
\begin{equation*}
\E|\hXn(u)|^2
=|\cTn|\qi \sum_{T\in\cTn} \biggl|\sum_{v\in T_n} e^{\ii \ell(v) u}\biggr|^2
\end{equation*}
and
\begin{equation}\label{samuel}
  \ftuu=\sum_{n=1}^\infty t^n|\cTn|\E|\hXn(u)|^2.
\end{equation}
Since 
$|\cTn|=[t^n]B(t)= \frac1{n+1}\binom{2n}{n} \sim \pi\qqi n^{-3/2} 4^n$,
\eqref{l3} is equivalent to
\begin{equation} \label{l3a}
[t^n] \ftuu \le \CC 4^n \frac{n\qq}{1+n u^4}, 
\qquad |u|\le\pi.
\end{equation}

We will prove this
using complex analysis.
We begin by studying $B$.

\begin{lemma}
  \label{LB}
$B=B(t)$ is a bounded analytic function of $t$ in the domain
  $\cD\=\bbC\setminus[1/4,+\infty)$.
Moreover, for $t\in\partial\cD=[1/4,+\infty)$, $B$ has continuous
boundary values $B_+(t)$ and $B_-(t)$ from the upper and lower side.
$B(t)$ (extended by $B_+$ or $B_-$)
is real if and only if $t\in(-\infty,1/4]$; on this interval
$B(t)$ is strictly increasing from $-1$ to $B(1/4)=1$.
\end{lemma}

\begin{proof}
  The first assertions are immediate from \eqref{B}.
Next, if $B(t)$ is real, then so is $t=B/(1+B)^2$.
It follows further from \eqref{B} that $B(t)$ is real for $t\le1/4$,
but $B_\pm(t)$ is not real for $t>1/4$.
The formula $t=B/(1+B)^2$ shows further that $B=-1$ is impossible, and
that $B=1$ if and only if $t=1/4$. Since $B(t)\to-1$ as $t\to-\infty$,
it follows by continuity that $-1<B(t)<1$ for $t<1/4$. For such $t$ we
have $dB/dt=(dt/dB)\qi=(1+B)^3/(1-B)>0$, which completes the proof.
\end{proof}

Let us for simplicity write $\ftu\=\ftuu$.

We first observe that, for any real $u$,
$\ftu$ is an analytic function of $t$ in the domain
$\cD'\=\cD\setminus(-\infty,-3/4]$.
Indeed, by \refC{CF2} and \refL{LB}, $\ftu$ is meromorphic in $\cD$
with poles when $1-B=0$ or $1+B-2B\cos u=0$.
In the first case, $B=1$ and thus $t=B(1+B)^{-2}=1/4$, which is
outside $\cD$.
In the second case, $B=1/(2\cos u-1)$. Since $2\cos u-1\in[-3,1]$,
this means that $B$ is real and either $B\ge1$ or $B\le-1/3$.
By \refL{LB}, $B\ge1$ is impossible in $\cD$, while $B\le-1/3$ implies 
$t=B(1+B)^{-2}\le-3/4$.

At this stage, we can apply, for any \emph{fixed} value of $u$,  the
standard results of {\em singularity
  analysis\/}~\cite{flajolet-odlyzko}. For $u=0$, we find
$$
[t^n] F_u(t) = \frac{4^n  n^{1/2}}{\sqrt \pi}\left(1+ O(1/n)\right),
$$
while for $u\neq 0$, 
$$
[t^n] F_u(t)= \frac 1{2(1-\cos u)^2} \frac{4^n  n^{-1/2}}{\sqrt
  \pi}\left(1+ O(1/n)\right).
$$
These results are certainly compatible with the desired
bound~\eqref{l3a}, but, as we need a \emph{uniform} bound, valid for
all $u$, we have to resort to the basic principles of singularity analysis.

By the Cauchy integral formula,
\begin{equation}
  \label{cauchy}
[t^n]\ftu = \frac1{2\pi\ii}\int_\Gamma \ftu \frac{\dd t}{t^{n+1}}
\end{equation}
for any contour $\Gamma$ in $\cD'$ that loops once around $0$.
We assume $n>4$ and choose 
a contour $\gG= \gG^{(n)}$ that depends on $n$: 
$\gG^{(n)}\=\gG_1\cup\gG_+\cup\gG_2\cup\gG_-$, where $\gG_1$ is
the circle $|t-1/4|=1/n$ in negative direction, $\gG_+$ and $\gG_-$ are
both the interval $[1/4+1/n,1/2]$, taken in opposite directions and
  using the boundary values $B_+$ and $B_-$, respectively,
and $\gG_2$ is the circle 
$|t|=1/2$
in positive direction.
(For convenience, we have pushed the contour to include part of the
boundary of $\cD'$; the reader who prefers staying strictly inside
$\cD'$ may replace $\gG_{\pm}$ by line segments close to the real
axis.)

Next we estimate $|\ftu|$ on $\gG^{(n)}$. 
\begin{lemma}\label{LF1}
For all $t\in\gG^{(n)}$ and $u\in[-\pi,\pi]$,
\begin{equation}\label{lf1}
  |\ftu| \le \CC \frac{n^{3/2}}{1+n u^4}.
\end{equation}
\end{lemma}

\begin{proof}
We claim that for $t\in\gG^{(n)}$
\begin{align}
  |1-B(t)|&\ge \cc n\qqi, \label{b1}
\\
  \bigl|1+B(t)-2\cos uB(t)\bigr|
&\ge \cc\max\bigpar{n\qqi,1-\cos u}. \label{b3}
\end{align}
The result then follows from \eqref{lf} and $1-\cos u\ge \cc u^2$.

In fact, \eqref{b1} is the special case 
$u=0$  
of \eqref{b3}, 
so it suffices to prove the latter.
Since, by compactness,  
$|B(t)|\ge\cc>0$ for $t\in \gG\subset\set{t:\frac1{20}\le|t|\le\frac12}$,
it is enough to prove
\begin{align}
  \bigl|B(t)\qi+1-2\cos u\bigr|
\ge \cc\max\bigpar{n\qqi,1-\cos u}. \label{b2}
\end{align}
Indeed,
\begin{equation}\label{bti}
  \begin{split}
  B(t)\qi+1-2\cos u
&
=\frac{1-2t+\sqrt{1-4t}}{2t}+1-2\cos u
\\&
=\frac{1-4t+\sqrt{1-4t}}{2t}+2(1-\cos u).	
  \end{split}
\end{equation}
For $t\in\gG_1$, this is $2\sqrtt+O(1/n)+2(1-\cos u)$.
Since $\Re\sqrtt\ge0$ and $1-\cos u\ge0$, then
\begin{equation*}
  \begin{split}
|B(t)\qi+1-2\cos u|
&\ge 
|2\sqrtt+2(1-\cos u)|-O(1/n)
\\&
\ge \max\set{2|\sqrtt|,2(1-\cos u)}-O(1/n)
\\&
= \max\set{4n\qqi,2(1-\cos u)}-O(1/n),
  \end{split}
\end{equation*}
which yields \eqref{b2}.

On $\gG_\pm$, $\sqrtt$ is imaginary, and 
\begin{equation}\label{bx}
\bigl|\Im\bigpar{  B(t)\qi+1-2\cos u}\bigr|
=\frac{|\sqrt{1-4t}|}{2t}
\ge
\sqrt{4t-1}
\ge 2 n\qqi.
\end{equation}
Moreover, if further $1-\cos u\le 2\sqrt{4t-1}$, 
\eqref{bx} also yields
$\bigl|\Im\bigpar{  B(t)\qi+1-2\cos u}\bigr|
\ge\sqrt{4t-1}
\ge\tfrac12(1-\cos u)$.
If, on the contrary, $1-\cos u> 2\sqrt{4t-1}$, then,
because $0\le4t-1\le1$ and thus $4t-1\le\sqrt{4t-1}$,
\begin{equation*}
  \begin{split}
\bigl|\Re\bigpar{  B(t)\qi+1-2\cos u}\bigr|
&
=2(1-\cos u)-\frac{4t-1}{2t}
\\&
\ge
2(1-\cos u)-2(4t-1)
\ge 1-\cos u.	
  \end{split}
\end{equation*}
In both cases, 
$\bigl|{ B(t)\qi+1-2\cos u}\bigr|
\ge\tfrac12(1-\cos u)$, which together with \eqref{bx} 
completes the verification of \eqref{b2} for $t\in\gG_\pm$.

Finally, for $t\in\gG_2$, we use compactness.
We observed above that $1+B(t)(1-2\cos u)=0$ is possible only for 
$t\in(-\infty,-3/4]\cup\set{1/4}$, and in particular not for $t\in\gG_2$; 
hence $\inf_{t\in\gG_2,\,u\in[0,2\pi]}|B(t)\qi+1-2\cos u|=\cc>0$, which
implies \eqref{b2} for $t\in\gG_2$.
\end{proof}

By \eqref{cauchy} and \eqref{lf1},
\begin{equation}\label{lf2}
[t^n]\ftu \le \int_\Gamma |\ftu|\, |t|^{-n-1}\,|d t|
 \le \CCx \frac{n^{3/2}}{1+n u^4}\int_\Gamma  |t|^{-n-1}\,|d t|.
\end{equation}
For $t\in\gG_1$, $|t|^{-n-1}=O(4^n)$, and thus 
$\int_{\Gamma_1}  |t|^{-n-1}\,|d t|=O(n\qi 4^n)$.
Secondly,
$\int_{\Gamma_\pm}  |t|^{-n-1}\,|d t| \le
\int_{1/4}^{1/2} t^{-n-1}\dd t 
\le n\qi 4^n.
$
Finally,
$
\int_{\Gamma_2}  |t|^{-n-1}\,|d t|=O(2^{n}).
$
Summing these estimates we find
$\int_{\Gamma}  |t|^{-n-1}\,|d t|=O(n\qi 4^n)$,
which together with \eqref{lf2} completes the proof of \eqref{l3a} and
thus \refL{L3} in the case of binary trees.

\smallskip
For complete binary trees, we use the well-known equivalence between 
binary and complete binary trees, where a binary tree $T$ of order $n$ is
identified with the internal nodes in a complete binary tree $T^c$ of
order $2n+1$. 
With this identification, one has
\begin{equation*}
X(j;T^c)= 
\begin{cases}
X(j-1;T)+X(j+1;T) & \hbox{if } j\neq 0,\\
1+X(-1;T)+X(1;T) & \hbox{otherwise.}
\end{cases}
\end{equation*}

Hence, temporarily using $X^c_n$ instead of
$X_n$ for the complete binary trees, it follows from \eqref{emma} that
$\widehat{X^c_{2n+1}}(u)=1+2\cos u \hXn(u)$.
Hence the estimate in \refL{L3} holds for complete binary trees too
(possibly with a different constant).

\section{Proof of Corollaries \ref{C2}--\ref{Cfise0}}\label{SpfC}

\begin{proof}[Proof of \refC{C2}]
This is immediate from \eqref{t2b} and the fact that $f_n\to f$ in
$\cor$ implies $f_n(j_n/n\qw)\to f(x)$, see \eg{} 
\cite[Theorem 5.5]{Billingsley}. 
\end{proof}

\begin{proof}[Proof of \refC{Cfise}]
By symmetry, $\fise(x)\eqd\fise(-x)$, so we may suppose $x\ge0$.
Then, as shown by Bousquet-M\'elou 
\cite[\S6.2.2, Conjecture 15 and Theorem 14]{BM},
for naturally embedded random binary trees,
$X_n(\floor{x n\qw})/n\qww\dto 2\qqi Y(2\qqi x)$ for a family
of random variables $Y(u)$, $u\ge0$, with \mgf{s}
$\E e^{aY(u)}=L(u,a)$. 
Combining this with \refC{C2}, we find 
$2\qwi\fise(2\qwi x)\eqd2\qqi Y(2\qqi x)$, $x\ge0$,
and thus
$\fise(x)\eqd2\qwi Y(2\qwi |x|)$.
(The normalization of $\fise$ in \cite{BM} is different.)
\end{proof}

\begin{proof}[Proof of \refC{Cfise0}]
By the proof of \refC{Cfise}, 
$\fise(0)\eqd2\qwi Y(0)$, where, by
\cite[Proposition 12 and Theorem 14]{BM}
$Y(0)\eqd 2\qq3\qi T\qqi$.

The (negative) moments of $T$ are given by the standard formula 
$\E T^{-s}=\Gamma(3s/2+1)/\Gamma(s+1)$, $s>-2/3$.
\end{proof}

\section{Other tree models}\label{Strees}

Consider a randomly labelled \cGWt{} as in \refConj{conj1}.
We know that the global limit result \eqref{sofie} holds, 
and the proof in \refS{Spf1} holds \emph{verbatim} in this case too
and shows that to prove \refConj{conj1}, it is sufficient to 
verify that the estimate of \refL{L3} holds.
We have not been able to do 
so in general, but we can show the required estimate
in the two special cases in \refT{Tplane}.

We consider thus in this section the two  families of 
\emph{labelled plane trees} that were studied 
 in \cite{BM}.   In the first family $\cTx1$, the root is labelled
 $0$, and the labels of 
two adjacent nodes differ by $\pm 1$. In the second family $\cTx2$, the latter
 condition is generalized by allowing the increments along edges to
 be $0, \pm 1$.

Again, we introduce a sequence of generating functions:
\begin{equation}\label{fkplane}
F_k(t,x_1,\dots,x_k) \equiv F(\bm x)
 \=
 \sum_{T\in \cT} t^{|T|} \prod_{i=1}^k \Bigpar{\sum_{v\in T} x_i^{\ell(v)}},
\end{equation}
where $\cT$ is either $\cTx1$ or $\cTx2$ and $|T|$ is the number of
{\em edges\/} in $T$.  
The following proposition is the counterpart, for each of the two new
families, of Proposition~\ref{propo-enum-binary}.

\begin{proposition}
\label{propo-enum-plane}
For plane trees with increments $\pm 1$,  the series $F_k$ can be
  determined by induction on 
  $k\ge 0$ using  
$$
F_k( \bm x)= 1 + t \sum_{I\subseteq [k]} \left(\prod_{i \in I} \bx_i
+\prod_{i \in I} x_i\right) F_{|I|} (x_I)F_{|J|} (x_J)
$$
where 
$J=[k] \setminus I$ and $\bx_i=1/x_i$.
For trees with increments $0, \pm 1$, the above equation becomes
$$
F_k( \bm x)= 1 + t \sum_{I\subseteq [k]} \left(1+\prod_{i \in I} \bx_i
+\prod_{i \in I} x_i\right) F_{|I|} (x_I)F_{|J|} (x_J),
$$
with the same notation as above. In both cases, each $F_k( \bm x)$
admits a rational expression in $F_0$ and the $x_i$. 
\end{proposition}
\begin{proof}
The proof is very similar to that of
Proposition~\ref{propo-enum-binary}. We now use the standard  recursive
description of  plane trees based on the deletion of the leftmost
subtree $T_1$ of a  tree $T$ (not reduced to a single node). This
leaves another plane 
tree $T_2$. Also, one has to take into account the
fact that  the label of the root of $T_1$ 
may now take two (or three) different values (depending on the family
of trees under consideration). Finally, the tree reduced to a single
node contributes $1$ in each $F_k$.  When $\cT=\cTx1$ and $k\ge 0$, this gives
\begin{multline*}
F_k(\bm x) = 1+ \sum_{T_1, T_2} t^{1+|T_1|+|T_2|} 
\prod_{i=1}^k \left( \sum_{v\in T_1} x_i^{\ell(v)-1}+ \sum_{v\in
  T_2} x_i^{\ell(v)} \right)\\
+ \sum_{T_1, T_2} t^{1+|T_1|+|T_2|} 
\prod_{i=1}^k \left(  \sum_{v\in T_1} x_i^{\ell(v)+1}+ \sum_{v\in
  T_2} x_i^{\ell(v)} \right),
\end{multline*}
and the result follows after expanding the products, and then
exchanging the sums.
\end{proof}

We easily find explicit formulas for $F_1$ and $F_2$ from 
\refP{propo-enum-plane}, \cf{} the proof of \refC{CF2}.
We leave the details to the reader and state only the result that we
need, in terms of the series $T= T(t)$ that counts  labelled trees
not reduced to a single node. Depending on which tree family is
studied, one has 
\begin{align*}
T=T^{(1)}&\= B(2t)=\frac{1-4t-\sqrt{1-8t}}{4t}  & &\text{for $\cTx1$},
\\
T=T^{(2)}&\=B(3t)=\frac{1-6t-\sqrt{1-12t}}{6t}  & &\text{for $\cTx2$}, 
\end{align*}
where the series $B(t)$ is defined by~\eqref{B}. 
\begin{corollary}
  \label{CF2plane}
For plane trees with increments $\pm 1$, 
\begin{equation}\label{lfp1}
\ftuu
=
\frac{(1+T)\left(1+T^2\cos^2u\right)}{(1-T)(1-T\cos u)^2}.
  \end{equation}
For plane trees with increments $0,\pm 1$,  
\begin{equation}\label{lfp01}
\ftuu
=
\frac{(1+T)(9+T^2\left(1+2\cos u)^2\right)}{(1-T)(3-T(1+2\cos u))^2}.
  \end{equation}
\end{corollary}

We may now complete the proof of \refT{Tplane} by 
the argument in \refS{Sgf}; we give a sketch only and leave again the
details to the reader.
First, 
the functions
$F_2(t/2,e^{\ii u},e^{-\ii u})$ (for $\cTx1$)
and $F_2(t/3,e^{\ii u},e^{-\ii u})$ (for $\cTx2$)
are analytic functions of $t\in\cD$ for every real $u$.
Next, in analogy with \refL{LF1},
with the same contour $\Gamma^{(n)}$ as there, 
for $t\in\Gamma$ and $|u|\le\pi$, 
\begin{align*} 
 \bigl|F_2(t/2,e^{\ii u},e^{-\ii u})\bigr| &\le \CC \frac{n^{3/2}}{1+n u^4}
& &\text{for $\cTx1$},
\\
 \bigl|F_2(t/3,e^{\ii u},e^{-\ii u})\bigr| &\le \CC \frac{n^{3/2}}{1+n u^4} 
&  &\text{for $\cTx2$}.
\end{align*}
Indeed, the proof is almost exactly the same; we replace the \lhs{} 
of \eqref{b3} by 
$\bigl|1-\cos uB(t)\bigr|$
and
$\bigl|3-B(t)-2\cos uB(t)\bigr|$ and similarly 
the \lhs{} \eqref{b2} by 
$\bigl|B(t)\qi-\cos u\bigr|$ and
$\bigl|3B(t)\qi-1-2\cos u\bigr|$, note the corresponding changes in
\eqref{bx} 
and argue as before.

The Cauchy integral formula \eqref{cauchy} then leads to, \cf{} \eqref{l3a},
for $|u|\le\pi$,
\begin{align*}
[t^n] \ftuu &\le \CC  8^n \frac{n\qq}{1+n u^4},  & &\text{for $\cTx1$},
\\
[t^n] \ftuu &\le \CC  12^n \frac{n\qq}{1+n u^4},  & &\text{for $\cTx2$}.
\end{align*}
By \eqref{samuel} and 
$|\cTx1_n|=2^n[t^n]B(t)\sim \pi\qqi n^{-3/2} 8^n$,
$|\cTx2_n|=3^n[t^n]B(t)\sim \pi\qqi n^{-3/2} 12^n$,
this yields \eqref{l3} for these two families.
(Note that we have let $|T|$ be the number of edges for $\cTx1$ and $\cTx2$;
thus we now should replace $X_n$ by $X_{n+1}$ in \eqref{samuel}, but
this makes no difference for \eqref{l3}.)

This completes the proof of \refT{Tplane}.
\label{Slastpf}

\section{Moments of the density of ISE}\label{Smom}

We know by \refC{Cfise} that $\fise(\gl)$ has a moment generating
function (defined in an interval containing 0), and thus finite
moments of all orders. We next present a formula for these moments, and more
generally for mixed moments involving several values of $\gl$.
We use a general method for occupation densities of Gaussian processes.
To state the formula, we introduce more notation.

Given $\zeta$ (which as always is $2\bex$), and $k$ points
$\ssk\in[0,1]$, the random vector $\bigpar{\tW(s_1),\dots,\tW(s_k)}$
has a Gaussian distribution with mean 0 and covariance matrix 
\begin{equation}
\gSzsk\=\bigpar{m(s_i,s_j;\zeta)}_{i,j=1}^k .  
\end{equation}
We let $\phizsk$ denote the density function of this distribution. 
(We may ignore the cases when the distribution is degenerate; \as{}
this happens only when $s_i=s_j$ for some $i$ and $j$.)

Using the construction in \refSS{SScrt} of the \BCRT{} from $\zeta$,
we can transfer these notations to the CRT.
Given $\zeta$ and $k$ points
$\xxk$ in the corresponding CRT, the random vector 
$\bigpar{\xW(x_1),\dots,\xW(x_k)}$
has a Gaussian distribution with mean 0 and covariance matrix 
\begin{equation}\label{gsx}
\gSzxk\=\bigpar{m(x_i,x_j;\zeta)}_{i,j=1}^k,
\end{equation}
where $m(x,y;\zeta)$ is the length of the common part of the paths 
from the root to $x$ and $y$ in the CRT.
We let $\phizxk$ denote the density function of this distribution, and
note that if $x_i=\tzeta(s_i)$, $i=1,\dots,k$, then
$m(x_i,x_j;\zeta)=m(s_i,s_j;\zeta)$ and $\phizxk=\phizsk$.

We further let $\XXk$ denote $k$ independent random nodes in the \BCRT.
\begin{theorem}
  \label{Tmom}
For any real numbers $\glk$,
\begin{equation}\label{tmom}
  \begin{split}
  \E\bigpar{\fise(\gl_1)\dotsm\fise(\gl_k)}
&=
\E
\intoi\dotsm\intoi \phizsk(\glk)\dssk
\\&=
\E \phizXk(\glk).
  \end{split}
\raisetag{\baselineskip}
\end{equation}
\end{theorem}

\begin{proof}
The equality of the last two expressions follows by the construction
of the \BCRT{} and the definitions above.

We define, for $\gl\in\bbR$ and $h>0$,
\begin{equation}\label{zh}
  \begin{split}
  Z_h(\gl)
&
=
h\qi\intoi\etta_{\bigsqpar{\tW(s)\in[\gl,\gl+h]}}\dd s
=
h\qi\muise[\gl,\gl+h]
\\&
=
h\qi\int_{\gl}^{\gl+h}\fise(y)\dd y.	
  \end{split}
\end{equation}
Since $\fise$ is continuous by \refT{T1}, $Z_h(\gl)\to\fise(\gl)$
\as{} as $h\to0$.
{}From this definition follows
\begin{align*}
 & \E\bigpar{Z_h(\gl_1)\dotsm Z_h(\gl_k)\mid\zeta}
\\
&
=
\intoi\dotsm\intoi h^{-k} 
 \P\bigpar{\tW(s_i)\in[\gl_i,\gl_i+h],\, i=1,\dots,k\mid\zeta}
\dssk
\\&
=
\intoi\dotsm\intoi h^{-k} 
\int_{\gl_1}^{\gl_1+h}\dotsm\int_{\gl_k}^{\gl_k+h} \phizsk(\yyk)\dyyk\dssk
\\&
=
\E\biggpar{h^{-k} 
\int_{\gl_1}^{\gl_1+h}\dotsm\int_{\gl_k}^{\gl_k+h} \phizXk(\yyk)\dyyk
\,\Big|\,\zeta}
\end{align*}
and thus
\begin{multline}\label{mom}
 \E\bigpar{Z_h(\gl_1)\dotsm Z_h(\gl_k)}
\\
=
\E h^{-k} 
\int_{\gl_1}^{\gl_1+h}\dotsm\int_{\gl_k}^{\gl_k+h} \phizXk(\yyk)\dyyk
.
\end{multline}

To obtain the conclusion, we now let $h\to0$; however, we have to
justify taking the limit inside the expectations on both sides.
For the right-hand side, we use the fact that a Gaussian distribution
in $\bbR^k$ with mean 0 has a density function that has its maximum at 0;
hence we can, by \refL{Lmom} below, use dominated convergence with
$\phizXk(0,\dots,0)$ as dominating function. Since $\phizXk$ is
continuous, the right-hand side of \eqref{mom} thus converges to the
right-hand side of \eqref{tmom}.

For the left-hand side we begin by applying Fatou's lemma, which now shows
that the left-hand side of \eqref{tmom} is \emph{at most} equal to the
right-hand side. By \refL{Lmom} below, this yields a uniform bound,
$C_k$ say, of the left-hand side for all $\glk$. It follows from
\eqref{zh} that 
$ \E\bigpar{Z_h(\gl_1)\dotsm Z_h(\gl_k)}\le C_k$ too, for every $h>0$.
If we here replace $k$ by $2k$, repeating every $\gl_i$ twice, we see
that the random variables $V_h\=Z_h(\gl_1)\dotsm Z_h(\gl_k)$ satisfy
$\E V_h^2\le C_{2k}$.
The variables $V_h$ are thus uniformly integrable, and from 
$V_h\to\fise(\gl_1)\dotsm\fise(\gl_k)$ as $h\to0$ follows 
$\E V_h\to\E\bigpar{\fise(\gl_1)\dotsm\fise(\gl_k)}$,
see \eg{} \cite[Theorems 5.4.2 and 5.5.2]{Gut}.
\end{proof}

\begin{lemma}
  \label{Lmom}
For every $k\ge1$,
$\E \phizXk(0,\dots,0)<\infty$.
\end{lemma}
\begin{proof}
  The subtree $R_k$ of the \BCRT{} spanned by $\XXk$ and the root $o$
  has $k-1$ internal nodes.
Let $R_k'$ be the subtree spanned by $o$ and the internal nodes of
  $R_k$, and let $\ell_1,\dots,\ell_k$ be the lengths of the $k$ edges
  that attach $\XXk$ to $R_k'$.
The values of $\xW$ along $R_k$ form a branching Brownian motion, \ie,
  $\xW$ is a Brownian motion along each edge of $R_k$ and all
  increments are independent. In particular, conditioned on $R_k$ and
  the values of $\xW$ on $R_k'$, the values $\xW(X_1),\dots,\xW(X_k)$
at the leaves are independent Gaussian variables with some means 
$b_1,\dots,b_k$ 
and  variances $\ell_1,\dots,\ell_k$.
The conditional density function is thus at most
  $\prod_{1}^k(2\pi\ell_i)\qqi$, 
and thus, taking the expectation and using \eqref{rk},
\begin{align*}
& \E \phizXk(0,\dots,0)\le \E \prod_{1}^k(2\pi\ell_i)\qqi
\\
&=(2k-3)!!\,(2\pi)^{-k/2}\!\idotsint \prod_{1}^k \ell_i\qqi 
\Bigpar{\sum_{i=1}^{2k-1}\ell_i} 
e^{-\tfrac12\bigpar{\sum_1^{2k-1}\ell_i}^2}
\dd\ell_1\dotsm\dd\ell_{2k-1}
\\
&<\infty. 
\end{align*}
\vskip-\baselineskip
\end{proof}

Since the distribution of the covariance matrix $\gSzXk$ is given by 
\eqref{gsx} and \eqref{rk}, it is in principle possible to write the \rhs{}
of \eqref{tmom} as a multiple integral. 
However, the expression becomes rather complicated for higher moments.
In the simplest case $\gl_1=\dots=\gl_k=0$,  \eqref{tmom} 
reduces to 
$\E \fise(0)^k=(2\pi)^{-k/2}\E \bigpar{\det\bigpar{\gSzXk}\qqi}$,
but even this seem difficult to compute in general. (These moments
were found by another method in \refC{Cfise0}.)

In the case $k=1$, \refT{Tmom} yields a simple formula 
for the average $\E\fise$ of the density, which equals 
the density of the average $\E\muise$, \ie{} 
the density of a random point chosen according to
the random ISE.
In the latter formulation, it was found by Aldous
\cite{AldousISE}.

\begin{corollary}
  \label{Cmom1}
For any real $\gl$,
\begin{equation*}
\E\fise(\gl)
=(2\pi)\qqi \intoo y\qq \exp\Bigpar{-\frac{\gl^2}{2y}-\frac{y^2}{2}}\dd y.
\end{equation*}
\end{corollary}
\begin{proof}
$\varphi_{\zeta;X_1}(\gl)=(2\pi y)\qqi e^{-\gl^2/(2y)}$, where
  $y=d(X_1,o)$, and $y$ has the density function $ye^{-y^2/2}$ by \eqref{rk}.
\end{proof}

Alternatively, expanding the Laplace transform of
Corollary~\ref{Cfise} in $a$ gives (see~\cite[Proposition 13]{BM}):
$$
\E\fise(\gl)
=\frac {2^{-1/4}} {\sqrt\pi}\sum_{m\ge 0}
  \frac{(-2^{3/4}|\lambda|)^m}{m!} \cos 
\frac{(m+1)\pi}4
\Gamma\left(\frac{m+3}4\right).
$$
Both expressions yield $\E\fise(0)=2\qwwi\pi\qqi\Gamma(3/4)$, as given by
\refC{Cfise0}. 

{}From \refC{Cmom1} follows easily by integration 
another formula by Aldous
\cite{AldousISE};
we leave the proof to the reader.

\begin{corollary}\label{Cmom2}
For every real $a>-1$,
\begin{align*}
  \E \intoooo|x|^a\dd\muise(x)
&=
  \E \intoooo|x|^a\fise(x)\dd x
=
\frac{2^{3a/4}}{\sqrt\pi}
 \Gamma\Bigpar{\frac a2+\frac12}\Gamma\Bigpar{\frac a4+1}.  
\end{align*}
\nopf
\end{corollary}
This extends~\eqref{2kth}.
(To see the equivalence when $a=2k$, use the duplication formula for
the Gamma function twice.) 


\section{The grand-moments of the ISE: Proofs}
\label{Smom3}

\begin{proof}[Proof of \refT{thm-grandmoments}]
Let $f$ be a continuous function on $\bbR$. 
First, if $f$ is bounded, then $\mu\mapsto\int f\dd\mu$ is a
continuous functional on
the space of probability measures on $\bbR$, and since
$\mu_n\dto\muise$ in this space, see \eqref{sofie},
it follows that
\begin{equation}
  \label{ej}
\int f\dd\mu_n \dto \int f\dd \muise.
\end{equation}
We need to extend this to unbounded $f$. Thus, let $f_u$, for $u>0$,
be the function that is equal to $f$ on $[-u,u]$, and is constant on
$(-\infty,-u]$ and on $[u,\infty)$. Since $f_u$ is bounded, \eqref{ej}
  applies to $f_u$, \ie{} $\int f_u\dd\mu_n \dto \int f_u\dd \muise$
  for every $u>0$.
Moreover, let 
$\WW_n\=\sup\set{|x|:x\in\supp\mu_n}=(2n)\qwi\sup\set{|\ell(v)|:v\in T_n}$.
By Marckert \cite[Theorem 5]{Marckert}, $\WW_n\dto \WW$ for some random
variable $\WW$. 
Consequently,
\begin{equation*}
  \begin{split}
  \limsup_\ntoo \P\Bigpar{\int f_u\dd\mu_n \neq \int f\dd\mu_n}
&
\le
  \limsup_\ntoo \P( \WW_n>u)
\le \P(\WW\ge u),	
  \end{split}
\end{equation*}
which tends to 0 as $u\to\infty$. Finally, 
$\int f_u\dd\muise \to \int f\dd\muise$ as $u\to\infty$, since
$\muise$ has compact support.
Consequently,
\cite[Theorem 4.2]{Billingsley} shows that \eqref{ej} holds for any
continuous $f$.

Taking $f(x)=x^k$ in \eqref{ej}, we obtain the convergence
$m_{k,n}\dto m_k$ of the moments asserted in
\refT{thm-grandmoments}. 
Moreover, taking $f$ to be a linear
combination of such monomials, we see that joint convergence holds
by the Cram\'er--Wold device \cite[Theorem 7.7]{Billingsley}.

In particular, for any partition $\gl$, $m_{\gl,n}\dto m_{\gl}$.
We will show that the expectation $\E( m_{\gl,n})$ converges as \ntoo.
Applying this to the partition $\gl'$ where each part in $\gl$ is
repeated twice, we see that also $\E (m_{\gl,n}^2)=\E (m_{\gl',n})$
converges. The variables $m_{\gl,n}$ are thus
uniformly integrable, and the limit of their expectations 
$\E (m_{\gl,n})$ equals the
expectation $\E( m_{\gl})$ of their limit,
see \eg{} \cite[Theorems 5.4.2 and 5.5.2]{Gut}.

To complete the proof of \refT{thm-grandmoments}, we thus have to show
that the grand-moments $\E (m_{\gl,n})$ of $\mu_n$ converge to the
limits stated in the theorem.
We introduce the
non-normalized moments of $\mu_n$:
$$
\bar M_{i,n} =  \sum_{v\in\tn}\ell(v)^i= 2^{i/4} n^{1+i/4}\, m_{i,n},
$$
as well as their factorial version, which is simpler to handle via
generating functions:
$$
M _{i,n} = \sum_{v\in\tn}\ell(v)\left(\ell(v)-1\right) \cdots
\left(\ell(v)-k+1\right).
$$
We also use the notation $M_{\lambda,n}$, analogous
to~\eqref{moments-notation}.
Then the remaining  part of \refT{thm-grandmoments} easily follows from:
\begin{proposition}\label{P12} 
As $n\rightarrow \infty$, the non-normalized factorial moments of
$\mu_n$ satisfy 
$$
\E( M_{\lambda, n}) =
\frac{\Gamma(1/2)n^{p+|\lambda|/4}
}{\Gamma(p+|\lambda|/4-1/2)} \left( c_\lambda  +o(1)\right)
$$
\end{proposition}
\begin{proof}
Let us first relate $M_{\lambda,n}$ to the generating functions of
Proposition~\ref{propo-enum-binary}. It is simple to see that
\begin{align}
\nonumber
\partial_\lambda F_p
&\= \frac{\partial^{|\lambda|}F_p}
{\partial x_1^{\lambda_1} \cdots \partial x_p^{\lambda_p}}(t,1, \
\ldots , 1)
=\sum_{T \in \cT} t^{|T|} \prod_{i=1}^p \left( 
\sum_{v \in T} \ell(v) (\ell(v)-1)\cdots(\ell(v)-\lambda_i+1)\right)
\\&
=\sum_{n\ge 0} t^n C_n \E(M_{\lambda,n}),\label{diff-lambda}
\end{align}
where $C_n=\binom{2n}{n}/(n+1)$ is the number of binary trees with
$n$ nodes, known as the $n$th Catalan number.
By Proposition~\ref{propo-enum-binary}, the series $\partial_\lambda
F_p$ is a rational function of $t$ and $\sqrt{1-4t}$. We want to study
the singularities of these series. We will prove
that, for $p>0$,
\beq
\label{Fp-form}
\partial_\lambda F_p= \frac{P_\lambda (t)+ Q_\lambda(t)
  \sqrt{1-4t}}{(1-4t)^{e_\lambda}}, 
\eeq
where $P_\lambda (t)$ and $Q_\lambda (t)$ are two Laurent polynomials
in $t$, 
and 
\begin{equation*}
e_\lambda= 
 p + \frac 1 2 \left\lfloor \frac {|\lambda|} 2 \right\rfloor -\frac 12
= 
\begin{cases}
\displaystyle 
p +  {|\lambda|}/ 4 - 1 /2 & \hbox{ if } |\lambda|
\hbox{ is even},\\
\displaystyle 
p +   {|\lambda|}/ 4 - 3/ 4 & \hbox{ if } |\lambda|
\hbox{ is odd.}
\end{cases}
\end{equation*}
(Note that $P_\gl$ and $Q_\gl$ may be singular at $t=0$, although
$\partial_\lambda F_p$ is analytic there.)
{}From~\eqref{Fp-form}, it follows that the only possible
singularity of $\partial_\lambda F_p$ is at $t=1/4$, and that, as
$t\rightarrow 1/4$,
$$
\partial_\lambda F_p =
\frac{c_\lambda+o(1)}{(1-4t)^{p+|\lambda|/4-1/2}} ,
$$
where 
$c_\gl=P(1/4)$ when $|\gl|$ is even and $c_\gl=0$ when $|\gl|$ is
odd. 
We will further show that
the numbers $c_\lambda$ satisfy the recurrence
relation~\eqref{c-lambda-rec}. The form~\eqref{Fp-form} and the above
singular behaviour do not hold when $p=0$, and should be replaced in
this case by the expression~\eqref{F0} of $F_0$ and the singular
behaviour 
$$
F_0=2-2\sqrt{1-4t} + O(1-4t).
$$

Assume for the moment that we have proved~\eqref{Fp-form}. Then 
the standard results of {singularity analysis}
\cite{flajolet-odlyzko} provide
$$
[t^n] \partial_\lambda F_p = C_n \E(M_{\lambda,n})
=
\frac{  4^n n^{p+|\lambda|/4-3/2}}
{\Gamma(p+|\lambda|/4-1/2)}\left(c_\lambda +o(1)\right).
$$
Given that $C_n\sim 4^nn^{-3/2}/\Gamma(1/2)$, this gives the result
stated in the proposition. Note that this asymptotic behaviour also 
holds for $p=0$, with $c_\emptyset =-2$.

Let us thus focus on~\eqref{Fp-form}. Our proof works by
induction on $p+|\lambda|$.

$\bullet$ If $p=0$, then $\lambda$ is the empty partition, and we have worked
out above the value of $F_0$ and its asymptotic behaviour when
$t\rightarrow 1/4$.

$\bullet$ If $p>0$ and $\lambda_1=0$, then 
$$\partial_\lambda F_p=\sum_{n\ge 0} t^n C_n \E(M_{\lambda,n})
= \sum_{n\ge 0} t^n C_n n \E(M_{\lambda',n})
= t \frac{\partial }{\partial t} \partial_{\lambda'}F_{p-1},
$$
where $\lambda'=(\lambda_2, \ldots , \lambda_p)$. Then the form~\eqref{Fp-form}
follows by a simple calculation from the induction hypothesis, and the
fact that 
$e_{\lambda}=e_{\lambda'} +1$. (We do not give the details.) This
calculation  also provides  the  value of
$c_\lambda$ in terms of $c_{\lambda'}$.  The case $p=1$ and $\lambda =
(0)$ has to be treated separately, since in that case
$\lambda'=\emptyset$ and the form~\eqref{Fp-form} is
not valid.

$\bullet$  If $p>0$ and $\lambda_1>0$, then all the parts of $\lambda$
are positive. Let us differentiate 
\eqref{fkx} $\lambda_1$ times with respect to
$x_1$, then $\lambda_2$ times with respect to $x_2$, and so on, and
then set $x_i=1$ in the result. Since $\lambda_i>0$ for all $i$, the
terms for which $I\cup J \neq [p]$ do not contribute, and we are
left with
\beq
\label{Fp-expr}
\partial_\lambda F_p = t \sum_{I\subseteq [p]}
\partial_{\lambda_I}\left(F_{|I|}(x_I)\prod_{i\in I} \bx _i \right)
\partial_{\lambda_J}\left( F_{|J|}(x_J) \prod_{j\in J} x _j\right)
\eeq
where $J=[p]\setminus I$ and, for any function $G(t,x_I)$, we denote
$$
\partial_{\lambda_I}G = \frac{\partial ^{|\lambda_I|} G}
{\partial x_{i_1} ^{\lambda_{i_1}} \cdots \partial x_{i_r}
  ^{\lambda_{i_r}}}(t,1,\ldots ,1)
$$
if  $I=\{i_1, \ldots , i_r\}$ with $i_1 <\ldots <i_r$. Now
$$
\partial_{\lambda_I}\left(F_{|I|}(x_I)\prod_{i\in I} \bx _i  \right)
=
\sum_{\sigma \le \lambda_I} (-1)^{|\lambda_I|-|\sigma|} \partial_\sigma
F_{|I|}
\prod_{i\in I}{\frac {\lambda_i !}  {\sigma_i!}
},
$$
where the sum runs over all non-negative $|I|$-tuples $\sigma=
(\sigma_i)_{i\in I}$ 
that are less than or equal to $\lambda_I$.
 The second derivative contains fewer terms:
$$
\partial_{\lambda_J}\left( F_{|J|}(x_J)\prod_{j\in J} x _j\right)
=
\sum_{\varepsilon } \partial
_{\lambda_J - \varepsilon} F_{|J|} \prod_{j\in J} \lambda_j^{\varepsilon_j}  ,
$$
where the sum runs over all $|J|$-tuples $\varepsilon= (\varepsilon_j)_{j\in
  J}$ such $\varepsilon_j \in\{0,1\}$ for all $j$.

Let us now bravely replace the two derivatives occurring
in~\eqref{Fp-expr} by their sum-expressions given above, and (mentally)
expand the product of these sums. This gives $\partial_\lambda F_p $
as a sum over $I$, $\sigma$ and $\varepsilon$. In this sum, the series
$\partial_\lambda F_p $ appears twice, namely 
\begin{itemize}
\item[(i)] for $I=\emptyset$,
$\sigma=\emptyset$ and $\varepsilon = (0, \ldots, 0)$, 
\item [(ii)]  for $I=[p]$, 
$\sigma=\lambda$ and $\varepsilon=\emptyset$. 
\end{itemize}
The corresponding summands
are the same in both cases, namely $tF_0\partial_\lambda F_p
$. Hence~\eqref{Fp-expr} can be rewritten as
$$
(1-2tF_0) \partial_\lambda F_p= t\sum_{I, \sigma, \varepsilon} {\rm SUMMAND},
$$
where the sum now excludes Cases (i) and (ii). In this sum, all
terms of the form $\partial_\tau F_k$ now satisfy $k+ |\tau| <p
+|\lambda|$, so that the induction hypothesis applies to them. Note
also that $(1-2tF_0)=\sqrt{1-4t}$, so that the previous equation
really reads
\beq\label{sqrtFp}
\sqrt{1-4t}\ \partial_\lambda F_p= t\sum_{I, \sigma, \varepsilon} {\rm
  SUMMAND}
={\rm RHS}.
\eeq
The latter observation is the key in our proof of~\eqref{Fp-form}.

In the right-hand side of the equation, let us study separately the
cases where $I$ or $J$ are empty.

\noindent{\bf First case:  $I$ or $J$ is empty.}
The contribution of the terms for which $I=\emptyset$ is 
\beq
\label{I-empty}
tF_0 \sum_{\varepsilon \neq 0}  \partial
_{\lambda - \varepsilon} F_{p} \prod_{j=1}^p \lambda_j^{\varepsilon_j} .
\eeq
The contribution of the terms for which $J=\emptyset$ (that is,
$I=[p]$) is
\beq
\label{J-empty}
tF_0 \sum_{\sigma < \lambda} (-1)^{|\lambda|-|\sigma|} \partial_\sigma
F_{p}
\prod_{i=1}^p \frac {\lambda_i !}  {\sigma_i!}.
\eeq
We observe that the  terms for which $|\varepsilon|=1$
in~\eqref{I-empty} cancel out with
the terms for which $|\sigma|=|\lambda|-1$ in~\eqref{J-empty}. (More
generally, the term 
associated with $\varepsilon$, when  $|\varepsilon|$ is odd, cancels out
with the term associated with $\sigma=\lambda-\varepsilon$, but we do not
need this property). After these cancellations, all the terms
$\partial_\tau F_k$ that appear in this part of ${\rm RHS}$ satisfy $k=p$
and $|\tau|\le |\lambda|-2$. In particular, $e_\tau \le e_\lambda
-1/2$ for each of them. The induction hypothesis then guarantees that
this part of ${\rm RHS}$ is of the form 
\begin{equation*}
{\rm RHS}_1=\frac{\PP1 (t)+ \QQ1(t)
  \sqrt{1-4t}}{(1-4t)^{e_\lambda -1/2}}, 
\end{equation*}
for two Laurent polynomials $\PP1(t)$ and $\QQ1(t)$. Given that we still
have to divide ${\rm RHS}$ by $\sqrt{1-4t}$ to obtain the expression of
$\partial_\lambda F_p$ (see~\eqref{sqrtFp}), this part of ${\rm RHS}$
is compatible with the expected form~\eqref{Fp-form}. 

Before turning our attention to the  case $\emptyset \neq I \neq
[p]$, let us work out the value of $\PP1(1/4)$, at least  when
$|\lambda|$ is even. In  ${\rm RHS}_1$, the only terms
$\partial_\tau F_p$  for which $e_\tau = e_\lambda
-1/2$ are those for which $|\tau|=|\lambda|-2$. That is, the terms for
which $|\varepsilon|=2$ in~\eqref{I-empty},
and the terms for which
$|\sigma|=|\lambda|-2$  
in~\eqref{J-empty}.
As $F_0\rightarrow 2$ when $t\rightarrow 1/4$, this means that
$$
{\rm RHS}_1 = \frac {1} {2(1-4t)^{e_\lambda -1/2}} 
\left( \sum_{1\le i<j\le p} c_{\lambda-\varepsilon_{i,j}}
\lambda_i \lambda _j 
+
 \sum _{\sigma \le \lambda, |\sigma|=|\lambda|-2}
c_\sigma \prod_{i=1}^p \frac{\lambda_i!}{\sigma_i!}+o(1) \right),
$$
where $\varepsilon_{i,j}$ is the $p$-tuple that has a one at positions $i$
and $j$, and zeros elsewhere.
In the second sum, the partitions $\sigma$ such that
$\sigma_i=\lambda_i-2$ for some $i$ contribute $\lambda_i
(\lambda_i-1)$, while those for which $\sigma_i=\lambda_i-1$ and
$\sigma_j=\lambda_j-1$, so that $\sigma= \lambda-\varepsilon_{i,j}$,
contribute $\lambda_i \lambda _j $, as in the 
first sum. A concise way of merging both sums consists in using the
notation of~\eqref{c-lambda-rec} and writing
$$
{\rm RHS}_1 = \frac {1} {(1-4t)^{e_\lambda -1/2}} \left(
\sum _{\sigma \le \lambda, |\sigma|=|\lambda|-2}
c_\sigma  \binom{\lambda}{\sigma} +o(1)\right),
$$
so that the polynomial $\PP1(t)$ satisfies
$$
\PP1(1/4)= \sum _{\sigma \le \lambda, |\sigma|=|\lambda|-2}
c_\sigma  \binom{\lambda}{\sigma},
$$ 
where we recognize the second part of~\eqref{c-lambda-rec}.

\medskip
\noindent
{\bf Second case: $I \neq \emptyset $ and $J \neq \emptyset $.}
In that case, the induction hypothesis~\eqref{Fp-form}  applies both
to  $\partial_\sigma
F_{|I|}$ and $
 \partial _{\lambda_J - \varepsilon} F_{|J|}$. Moreover,
\beq
\label{ineq-e}
e_\sigma+ e_{\lambda_J - \varepsilon} \le e_{\lambda_I} +  e_{\lambda_J}
= p + \frac 1 2 \left\lfloor \frac{|\lambda_I|} 2\right\rfloor
+ \frac 1 2 \left\lfloor \frac{|\lambda_J|} 2 \right\rfloor -1
\le e_\lambda -\frac 1 2.
\eeq
This implies that the part of ${\rm RHS}$ for which 
$\emptyset\neq I \neq [p]$ can be written as
$$
{\rm RHS}_2= 
\frac{\PP2 (t)+ \QQ2(t)
  \sqrt{1-4t}}{(1-4t)^{e_\lambda -1/2}}, 
$$
for two Laurent polynomials $\PP2(t)$ and $\QQ2(t)$. Given that  ${\rm RHS}_1$
has also this form, we can conclude at last that~\eqref{Fp-form} holds. 

 Let us finally work out the value of $\PP2(1/4)$, at least  when
$|\lambda|$ is even. The only way for the inequalities~\eqref{ineq-e}
 to be equalities is to take $\sigma=\lambda_I$, $\varepsilon
 =0$, with
 $|\lambda_I|$ and $|\lambda_J|$ even. Going back to~\eqref{Fp-expr},
 this means that the dominant contribution in ${\rm RHS}_2$ is given by
$$
\frac1{4 (1-4t)^{e_\lambda -1/2}} \sum_{\emptyset \neq I \subsetneq [p]}
  c_{\lambda_I}c_{\lambda_J}.
$$
In other words, 
$$\PP2(1/4)=\frac 1 4  \sum_{\emptyset \neq I \subsetneq [p]}
  c_{\lambda_I}c_{\lambda_J},
$$
which gives the first part in~\eqref{c-lambda-rec}. 
This completes the proof of \refP{P12} and thus of \refT{thm-grandmoments}.
\end{proof}
\noqed
\end{proof}

The proof of \refT{thm-grandmoments-plane} is almost the same, using 
the generating functions and recursion relations in
\refP{propo-enum-plane} and replacing $1-4t$ by $1-8t$ and $1-12t$,
respectively.
We omit the details.

To conclude this section, let us sketch the proof of
Theorem~\ref{thm-grandmoments-exc}. 
We define the generating functions $F_k$ as in \eqref{fk}, but
replacing $\ell(v)$ by the depth $d(v)$.
Then \eqref{fkx} holds with $\bx_i$ replaced by $x_i$.
In particular, $F_0$ is still given by~\eqref{F0} 
and each $F_k( \bm x)$ admits a rational expression in $F_0$ and the $x_i$.
We claim that then \eqref{Fp-form} holds, with 
$e_\gl=p+|\lambda|/2-1/2$, and 
$$
\partial_\lambda F_p =
\frac{d_\lambda+o(1)}{(1-4t)^{p+|\lambda|/2-1/2}}.
$$
This is proved by induction as above. Note that, after
\eqref{Fp-expr},
the $\partial_{\gl_I}$ term expands exactly as the $\partial_{\gl_J}$
term,
without $(-1)^{|\gl_I|-|\sigma|}$  and thus without cancellation;
this ultimately explains why the exponents $e_\lambda$ increase faster
for the horizontal profile than for the vertical.
The rest is as above.

The same applies to plane trees, with $F_k$ defined as in 
\eqref{fkplane} with $\ell(v)$ replaced by $d(v)$, and the recursion relation
$$
F_k( \bm x)
 = 1 + t \sum_{I\subseteq[k]} \left(\prod_{i \in I} x_i\right)
 F_{|I|} (x_I)F_{|J|} (x_J),
$$
where $J=[k]\setminus I$. We omit the details.



\newcommand\AAP{\emph{Adv. Appl. Probab.} }
\newcommand\JAP{\emph{J. Appl. Probab.} }
\newcommand\JAMS{\emph{J. \AMS} }
\newcommand\MAMS{\emph{Memoirs \AMS} }
\newcommand\PAMS{\emph{Proc. \AMS} }
\newcommand\TAMS{\emph{Trans. \AMS} }
\newcommand\AnnMS{\emph{Ann. Math. Statist.} }
\newcommand\AnnPr{\emph{Ann. Probab.} }
\newcommand\CPC{\emph{Combin. Probab. Comput.} }
\newcommand\JMAA{\emph{J. Math. Anal. Appl.} }
\newcommand\RSA{\emph{Random Struct. Alg.} }
\newcommand\ZW{\emph{Z. Wahrsch. Verw. Gebiete} }
\newcommand\DMTCS{\jour{Discr. Math. Theor. Comput. Sci.} }

\newcommand\AMS{Amer. Math. Soc.}
\newcommand\Springer{Springer}
\newcommand\Wiley{Wiley}

\newcommand\vol{\textbf}
\newcommand\jour{\emph}
\newcommand\book{\emph}
\newcommand\inbook{\emph}
\def\no#1#2,{\unskip#2, no.~#1,} 

\newcommand\webcite[1]{\hfil\penalty0\texttt{\def~{\~{}}#1}\hfill\hfill}
\newcommand\webcitesvante{\webcite{http://www.math.uu.se/\~{}svante/papers/}}
\newcommand\arxiv[1]{arXiv:#1}

\def\nobibitem#1\par{}

\end{document}